\documentclass[11pt]{amsart}

\usepackage{amsmath}
\usepackage{amssymb}
\usepackage{latexsym}
\usepackage{bm}

\newtheorem{theo}{Theorem}[section]
\newtheorem{prop}[theo]{Proposition}
\newtheorem{lem}[theo]{Lemma}
\newtheorem{rem}[theo]{Remark}

\newtheorem{cor}[theo]{Corollary}

\newcommand{\R}{{\Bbb R}}
\newcommand{\N}{{\Bbb N}}

\newcommand{\T}{{\Bbb T}}

\renewcommand{\S}{{\Bbb S}}

\newcommand{\Sup}{\displaystyle \sup}

\newcommand{\E}{{\Bbb E}}

\newcommand{\M}{{\mathcal M}}
\renewcommand{\L}{{\mathcal L}}
\renewcommand{\o}{\omega}
\renewcommand{\H}{{\mathcal H}}

\newcommand{\U}{{\mathcal U}}
\def\eps{\varepsilon}

\def\b{\beta}
\def\be{\beta}\def\de{\delta}
\def\De{{\bm \delta}}
\def\ka{\kappa}
\def\th{\theta}
\def\Th{\Theta}
\def\pp{\frac{\tilde \psi^n}{\psi_\infty}}
\def\a{\alpha}

\renewcommand{\div}{{\rm div}}

\def\les{\lesssim}

\def\o{\overline}

\begin{document}

\title{ Global existence of weak solutions to the 
  FENE dumbbell model of polymeric flows  }
\author
{
Nader Masmoudi} 
\address{
 Courant Institute, New York University \\
  251 Mercer St,
New York NY 10012\\
email:masmoudi@cims.nyu.edu }

\date{}

\maketitle



\vspace{8mm}

\noindent \underline{\bf Abstract}

Systems coupling fluids and  polymers 
 are of great interest in many branches of sciences. 
One of the  models  to describe them is  the FENE (Finite Extensible Nonlinear Elastic)
  dumbbell model.
 We prove global existence of weak solutions to  the FENE 
 dumbbell model of polymeric flows for a very general class of 
potentials. The main problem  is the passage to the limit in 
 a nonlinear term that has no obvious compactness properties. 
 The proof uses many weak convergence techniques. 
In particular it is based on the control of the propagation of strong 
convergence  of some well chosen quantity by studying  a transport equation for 
its    defect measure.   


\section{introduction}

Systems coupling fluids and  polymers 
 are of great interest in many branches of applied physics,  chemistry
and biology. They are of course used in many industrial and medical 
 applications
such as food processing, blood flows... 
 Although a  polymer molecule may be a very complicated object, there are 
simple theories to model it. 
One of these models  is the FENE (Finite Extensible Nonlinear Elastic)
  dumbbell model. In this model, a polymer 
is idealized as an ``elastic dumbbell''  
consisting of two ``beads'' joined by a spring which can be 
represented by a vector $R$ (see Bird, 
		   Curtis,  Amstrong and  Hassager
 \cite{BAH77,BCAH87}, 
Doi  and Edwards \cite{DE86} for some physical introduction to the
model
   and  Ottinger \cite{Ottinger96} for a more mathematical treatment 
(in particular the stochastic point of view)
   of it and Owens and Phillips \cite{OP02} for the  computational aspect).   
In the FENE model \eqref{micro}, the polymer elongation  $R$  cannot exceed a limit 
$R_0$. This yields some nice mathematical problems near the boundary,
namely 
when $|R|$ approaches $R_0$.   
At the level 
of the polymeric liquid, we get a system coupling the Navier-Stokes 
equation for the fluid velocity with a Fokker-Planck equation 
describing the evolution of the polymer density. This density depends 
on $t,x$ and $R$.  The coupling 
comes from an  extra stress term in the fluid equation due 
to the microscopic effect of the  polymers. This is the micro-macro
interaction. There is also a  drift term in the Fokker-Planck equation  
that depends on the spatial gradient of  the velocity. 
This is a  macro-micro term. The coupling 
satisfies the fact that the free-energy 
dissipates which is important from the physical point of view.  Mathematically, 
this is also  important to get uniform bounds and hence prove global existence 
of weak solutions. 

 The system obtained  attempts to describe the behavior of this 
complex mixture of polymers  and fluid, and as such, it  presents 
 numerous challenges, simultaneously at the level of 
their derivation \cite{DLP02}, the level of their numerical simulation  
\cite{OP02,Keunings89}, the level of their physical properties (rheology) and that of their
 mathematical  treatment (see references below).
  In this paper we concentrate on 
the mathematical treatment and more precisely the global existence 
of weak solutions to the 
  FENE dumbbell model (\ref{micro}).  These solutions  are the generalization 
 of the Leray weak solutions \cite{Leray34am,Leray33} of the incompressible Navier-Stokes system 
to the FENE model.

An approximate closure of the linear Fokker-Planck equation reduces the 
description to a  closed viscoelastic equation for the 
added stresses themselves. This leads to well-known  non-Newtonian fluid models
such as the Oldroyd B model or the FENE-P model (see for instance 
\cite{DLY05,DLP02}).   
 These models 
   have been studied extensively.  Guillop{\'e}  and Saut \cite{GS90, GS90b} proved 
the existence of local strong solutions,   
Fern{\'a}ndez-Cara,  Guill{\'e}n  and
             Ortega  \cite{FGO98},   \cite{FGO97} and \cite{FGOb}
proved local well posedness in Sobolev spaces. In Chemin and Masmoudi 
\cite{CM01} local and global well-posedness in critical Besov spaces 
was given. 
For global existence of weak solutions, 
 we refer to Lions and Masmoudi  \cite{LM00cam}. 
We also  mention Lin, Liu and  Zhang  \cite{LLZ05} where 
a formulation based on the deformation tensor is used to 
study the Oldroyd-B model. Global existence for small data 
was also proved  in \cite{LZ05,LLZ08}.

At the micro-macro level, there are also several works. Indeed,  
from the  mathematical point of view, the FENE model and  some simplifications 
of it  were  studied by several authors.
In particular Renardy \cite{Renardy91} proved the local existence 
in Sobolev space where the potential $\U$ is given by $\U (R) = (1-|R|^2)^{1-\sigma}$ for 
some $\sigma > 1$. 
 W. E, Li and Zhang \cite{ELZ04} proved local existence when $R$ is taken in the 
whole space and under some growth condition on the potential. Also,   
Jourdain, Lelievre and Le Bris \cite{JLL04} proved local existence 
in the case $b=2k>6$   for  a Couette flow by solving a stochastic differential equation
(see also \cite{JLLO06} for the use of entropy inequality methods  to prove 
exponential convergence to equilibrium).    
   Zhang and Zhang \cite{ZZ06} 
 proved local well-posedness for the FENE model when $b>76$. Local well-posedness 
was also proved in \cite{Masmoudi08cpam} when $b=2k > 0$ (see also \cite{KP10}). 
One of the main 
ingredients of   \cite{Masmoudi08cpam}  is the use of Hardy type
inequalities  to control 
the extra stress tensor by the $H^1$ norm in $R$ which comes from the
diffusion in $R$. In particular no regularity in $R$ is necessary for
the 
initial data.  
 Moreover,  Lin, Liu and  Zhang  \cite{LLZ07}
  proved global existence near equilibrium  under some 
restrictions on the potential (see also the related work \cite{LLZ08}). 
 Recently many other works dealt with different aspect of the system. 
In particular the problem in a thin film was considered in \cite{Chupin09}, 
the problem of the long time behavior was considered in \cite{Schonbek09,LP08,ACM09}, 
the problem of global existence in smooth spaces  in 2D for some simplified 
models (when there is a bound on $\tau $ in $L^\infty$) was considered in 
\cite{CFTZ07,LZZ08,CM08,MZZ08}, the problem of non-blow up criterion was considered in 
 \cite{LMZ10}, the problem of stationary solution was considered in 
\cite{Chupin09,Chupin09prep}, the study of the boundary condition at $\partial B$ was 
considered in \cite{LS10,LL08}.

 More related to this paper, the  construction of  global weak solutions
for simplified models was considered in \cite{BSS05,BS07,BS08,ZZZ08,Schonbek09,BS10prep}  in 
the case the system  is regularized by some diffusion in the space
variable  or by  a microscopic cut-off. 
The case of the  co-rotational  model was considered in \cite{LM07}. 
The  co-rotational  model preserves
some of the compactness difficulties of the full model. It allows to
get more integrability on the $\psi$ which makes the compactness
analysis 
much simpler.

We end this introduction by mentioning other micro-macro models. Indeed, 
  a principle based on an energy dissipation balance was proposed in \cite{Constantin05}, 
where the regularity of  nonlinear Fokker-Planck systems coupled with Stokes equations 
in 3D was also proved. In particular the Doi model (or Rigid model) was considered 
in \cite{otto} where   the linear
Fokker-Planck system is  coupled with a stationary  Stokes equations.   
The nonlinear Fokker-Planck equation driven by a time averaged 
Navier-Stokes system in 2D was studied in \cite{CFTZ07} (see also \cite{CM08}). 
Recently, there were many review papers dealing with different mathematical aspects 
of these models \cite{Renardy00b,LZ07,LL09}. In particular we refer to \cite{LL09}
for an   exhaustive list of references dealing with the numerical point of view.

\subsection{The FENE model}
A  macro-molecule is idealized as an ``elastic dumbbell''  
consisting of two ``beads'' joined by a spring which can be 
modeled by a vector $R$ (see \cite{BCAH87}). Before writing 
our main system \eqref{micro}, let us discuss the main physical 
assumptions that lead to it: 

\begin{itemize}
\item The polymers are described by their density at each time t, position x 
and elongation $R$.  This is a  {\it kinetic description}  of the polymers. 

\item The inertia of the polymers is neglected and hence the sum of the forces  
applied on each polymer vanishes. We refer to \cite{DL09} where inertia is taken 
into account. Moreover,  the limit $m$ goes to zero it studied where $m$ is 
the mass of the beads. 

\item The polymer solution is supposed to be dilute and hence the interaction 
between different polymers is neglected. This is why we get a linear 
 Fokker-Planck equation. Let us also mention that there are models for 
polymer  melts such as the reptation model (see for instance \cite{Ottinger96}).

\item The polymer is described by one vector $R$ in $B(0,R_0)$.
 Let us mention that there 
are models where each polymers is described by one vector $R$  such that 
$|R|=1$ (the rigid case, see \cite{CM08}) or by $K$ vectors $R_i$, $1\leq i \leq K$ (see
\cite{BS10prep}). Usually the
difference between these models comes from the length of the polymers as 
well as their electric properties. 

\item In the Fokker-Planck equation an upper-convected derivative is used. 
This is can be seen as the most physical one. Other used derivatives are 
the lower-convected and the co-rotational ones (see \cite{BAH77,BCAH87}). The co-rotational one 
has the mathematical advantage that one has better a priori estimates (see \cite{LM07}).

\item  We neglect the diffusion in $x$ in the   Fokker-Planck  equation. Indeed, 
this diffusion is much smaller than the diffusion in $R$. Actually, it makes  
the mathematical problem much simpler.

\end{itemize}

Under these assumptions, the micro-macro approach consists in  writing a coupled 
multi-scale system :  
\begin{equation} \label{micro}   \left\{  \begin{array}{l}
  {\partial_t u} + (u\cdot \nabla) u- \nu \Delta u + \nabla p
 = {\div} \tau, \quad {\div} u = 0, \\
\\
\partial_t \psi + u. \nabla \psi =   {\rm div}_R \Big[ - \nabla u  \,  R \psi
      + {\beta} \nabla \psi +   \nabla \U  \psi  \Big]  \\
\\
\tau_{ij} =   \int_B   (R_i \otimes  \nabla_j \U) 
 \psi(t,x,R) dR \,  \quad \quad 
 ( \nabla \U  \psi  +  {\beta}  \nabla \psi).n = 0  \;  \hbox{on} \;   
 \partial B(0,R_0).       
  \end{array} \right. 
\end{equation}   
 
In (\ref{micro}), $\psi(t,x,R)$   denotes the distribution function for the 
internal configuration and $F(R) = \nabla_R \U$ is the spring  force 
 which derives from a potential $\U$ and $\U(R) = - {k }   {\rm log} (1 - |R|^2 /|R_0|^2  )$ for some 
constant $k> 0 $. Besides, $\beta$ is related to the 
temperature of the system and $\nu>0$ is the viscosity of the fluid.
 In the sequel, we  will take  $\beta = 1$.

Here,  $R$ is  in a bounded  ball $B(0,R_0)$ of radius $R_0$ which means that 
the extensibility of the polymers is finite and $x \in \Omega$  where $\Omega$ is a bounded 
domain of $\R^D$ where $D \geq 2$ or $\Omega = \T^D$ or $\Omega = \R^D$. In the case 
$\Omega$ has a boundary, we add the Dirichlet boundary condition $u=0$ on $\partial \Omega$.  
 We have also    
 to add a boundary condition to insure the conservation of $\psi$, namely 
$( -\nabla u R \psi +   \nabla \U  \psi  +  {\beta } \nabla \psi).n = 0 $ on $\partial B(0,R_0)$. 
The boundary condition on $\partial B(0,R_0)$ insures the conservation 
of the polymer density and
should be understood in the weak sense, namely for any function $g(R) \in C^1(B)$, 
we have 
\begin{equation} 
\partial_t \int_B g \psi dR  + u. \nabla_x \int_B g \psi dR  =  
   - \int_B \nabla_R  g \Big[ -  \nabla u  \,  R\,  \psi
      + {\beta} \nabla \psi +   \nabla \U  \psi  \Big] dR.  
\end{equation}
Notice in particular that it implies that $\psi = 0 $ on $\partial B(0,R_0)$
and that  if initially  $\int \psi(t=0,x,R) dR = 1$, then for all $t$ and $x$, we have 
$\int \psi(t,x,R) dR = 1$.
We will see later an other way of understanding this singular 
boundary condition.

When doing numerical simulation on the FENE model, it is usually better 
to think of the distribution function $\psi$ as  the density of a random variable 
$R$ which solves (see \cite{Ottinger96})
\begin{equation} \label{stoch}
dR + u. \nabla R dt = ( \nabla u R - \nabla_R \U (R)   )dt  + \sqrt2 dW_t  
\end{equation}
where the stochastic process $W_t$ is the standard Brownian motion in $\R^N$ 
and the additional stress tensor is given by the following expectation 
$\tau = \E(R_i \otimes  \nabla_j \U) $. 
Of course, we may need to add a boundary condition for (\ref{stoch}) if 
$R$ reaches the boundary of $B$. This is done by requiring that 
$R$ stays in $\overline B$ (see \cite{JL03}).  
Using this stochastic formulation  has the advantage of replacing the 
second equation of (\ref{fene}) which has $2D+1$ variables by (\ref{stoch}). 
Of course one has to solve (\ref{stoch}) several times to get the 
expectation $\tau$ which is the only information needed in the 
fluid equation. This strategy was used for instance by 
Keunings  \cite{Keunings97} (see also \cite{GHLKL99}) and by 
 {\"O}ttinger \cite{Ottinger96} (see also \cite{GO97}).

In the sequel, we will only deal with the FENE model and 
 we will take $\beta =1$ and $R_0=1$.

\section{Statement of the results}

This paper is devoted to the proof of global existence of 
free-energy weak solutions  to the FENE model. 
The main difficulty of the construction 
 is the passage to the limit in an approximate system   in the 
nonlinear term $\nabla u^n \psi^n$. Indeed, we only have a
uniform bound on $\nabla u^n$  in $L^2((0,T)\times \Omega)$ 
and $\psi^n$ in $L^\infty((0,T)\times \Omega; L^1(B))$ for all $T>0$ and so
assuming that $u^n$ and $\psi^n$ converge weakly to $u$ and $\psi$,  
it is not clear how to deduce that  $\nabla u^n \psi^n$ converges 
weakly to $\nabla u \psi $.

Before mentioning our main result, let us recall that 
  the  construction of  global weak solutions
to  simplified models was considered in \cite{BS07,BS08,Schonbek09,LM07,ZZZ08}. 
In particular in \cite{BS07}  a diffusion  in the space variable 
in the $\psi$ equation is added.  Mathematically this yields 
a bound on $\nabla_x \sqrt{\psi}$ in  $L^2((0,T)\times \Omega \times B)$
and hence one can easily pass to the limit in the product $\nabla u^n \psi^n$
using the Lions-Aubin lemma.  
  This extra diffusion  term is physically justifiable 
but it is much smaller than  the diffusion in the $R$ variable and this 
is why we did not include  it here. Recently, Barrett and Suli
\cite{BS10prep} extended 
their results to the case of bead-spring chain models where each polymer is 
described by $K$ springs $R^i$, $1\leq i \leq K$  
 again with diffusion in the $x$ variable.   
Also, in \cite{LM07}, the co-rotational model was considered. It allowed us 
to get more a priori estimates on $\psi^n$, namely one can get 
that $\psi^n$ is in all $L^p$ spaces. An argument based on propagation 
of compactness similar to the one used in \cite{LM00cam} allowed us to 
conclude.   
 
Here, we consider the more physical model \eqref{micro}.  
The system (\ref{micro}) has to be complemented with an initial data 
$u(t=0) = u_0$ and $\psi(t=0) = \psi_0$. 

 
Notice that $(u=0, \psi_\infty)$ where $\psi_\infty $
\begin{equation}
\psi_\infty(R)  = \frac{e^{-\U(R)}}{ \int_B e^{-\U(R')} dR'} 
 \end{equation}
defines a stationary solution of (\ref{micro}). To state our result,  we
first impose some conditions on the initial data. We take 
 $u_0 (x)   \in  L^2  (\Omega) $, div$(u_0) = 0  $ 
 and $\psi_0 (x,R) \geq 0  $ such that 
$\rho_0(x) =\int \psi_0 dR  \in L^\infty(\Omega)  $. Here $\rho_0(x)$ is 
the initial density of polymers at the position $x$.  We also assume 
the following entropy bound :   
$ \frac{\psi_0}{\rho_0 \psi_\infty}  \in  L\log L  (\Omega
\times B ,   dx  { \rho_0(x)  \psi_{\infty}  dR}  )  )  $, namely 
\begin{equation}  \label{ini-ent}
\|  \frac{\psi_0}{\rho_0 \psi_\infty}   \|_{
 L\log L  (\Omega
\times B ,     {\rho_0(x)   \psi_{\infty}   dR dx}  ) }  
=  
\int\int_{\Omega
\times B }   (   \frac{\psi_0}{\rho_0 \psi_\infty}    \log
\frac{\psi_0}{\rho_0 \psi_\infty}    - 
   \frac{\psi_0}{\rho_0 \psi_\infty}   + 1 ) \rho_0(x)   \psi_{\infty}  dR dx  <
\infty.  
\end{equation}
Finally, we also assume the following   $L^{1/2}_x   L\log^2 L  $  bound, that we will call 
   ``$\log^2$''  bound: 
\begin{equation}  \label{ini-log2}
 \int_\Omega  \frac{\int_B \psi_0  \log^2 \frac{\psi_0}{\rho_0
    \psi_\infty}    }  { 1 + \left[  \int_B \psi_0  \log^2 \frac{\psi_0}{\rho_0
    \psi_\infty}    \right]^{1/2} }  dx < \infty.  
\end{equation}
Notice that interpolating \eqref{ini-log2} with the $L^\infty$ bound on $\rho_0$, we
can deduce the $ L\log L  $ bound \eqref{ini-ent}. 

\begin{theo} 
\label{fene}  
Take a divergence free field 
 $u_0 (x)   \in  L^2  (\Omega) $ and $\psi_0 (x,R) \geq 0  $ such that 
$\rho_0(x) =\int \psi_0 dR  \in L^\infty(\Omega)  $ and
\eqref{ini-ent} and \eqref{ini-log2} hold. 
Then, \eqref{micro}   has a global weak solution  
  $(u,\psi)$   such that  $ u \in  L^\infty(\R_+; L^2) \cap L^2 (\R_+; \dot  H^1) 
 $,  $  \frac{\psi}{\rho \psi_\infty} \in 
  L^\infty( \R_+;  L\log L  (\Omega
\times B ,   dx  { \rho(x)   \psi_{\infty}  dR}  )  )    ) $ 
where  $\rho(x) = \int_B \psi dR $   and    $ \sqrt{ \frac{\psi}{\psi_\infty}   } \in 
  L^2( \R_+; L^2 (\Omega;   \dot  H^1_R  ( { \psi_{\infty}} dR   )  ) ) $
  and
  \eqref{free} holds with an inequality $\leq$ instead of the equality
  and  \eqref{log2-b} holds (with $\Omega$ replaced by any compact $K$
  of $\Omega$ in the whole space case).

\end{theo}

\begin{rem}
1) Of course $u$ and $\psi$ have also some time regularity in some 
negative Sobolev spaces in $x$ and $R$.  This allows to give a sense to the
initial data (see \cite{LM00cam} for more details). 

2) By   $f  \in L\log L  (\Omega
\times B ,   dx  {\rho(x)   \psi_{\infty}  dR}  )   $ we  mean  that 
$\int\int_{\Omega
\times B }   (f  \log f  - f + 1 ) \rho(x)   \psi_{\infty}  dR <
\infty  $. Notice that \eqref{ini-ent} does not  really define a norm. One
can of course define a norm using Orlicz spaces. However, we do not
need to do it here.  

3) If the domain $\Omega$ has finite measure (bounded domain  or torus)  then, 
 the extra bound \eqref{ini-log2} reduces to 
  $\int_\Omega \left[  \int_B \psi_0  \log^2
  \frac{\psi_0}{\rho_0   \psi_\infty}    \right]^{1/2} dx  < \infty. $
This extra bound  on the
initial data allows us to prove the extra  bound \eqref{log2-b} 
 on the solution. This is useful to get some sort of
 equi-integrability of the extra stress tensor. 
Of course this is a very mild extra assumption, but it would be 
  nice to see if one can prove the same result without it. Moreover,
due to the local character of the weak compactness proof,   the 
assumption \eqref{log2-b}  can be weakened by assuming the bound to
hold
locally in space, namely  $\int_K  \left[  \int_B \psi_0  \log^2
  \frac{\psi_0}{\rho_0   \psi_\infty}    \right]^{1/2} dx  < \infty $
for any compact set  $K$ of $\Omega$. 

4) For the simplicity of the presentation, the proof will be given 
in the case $\rho_0(x) $ is constant equal $1$ and $\Omega$ has finite
measure. We will also indicate 
the necessary changes to be done in the general case.


\end{rem}


The paper is organized as follows. In the next section, we give 
some preliminaries where we prove  some Hardy type inequalities. 
 In section  
\ref{aprio}, we derive  some a priori estimates for the full model 
(\ref{micro}). In particular we recall the free energy estimate 
as well as a new ``$log^2$''    a priori estimate  which is 
useful in controlling the transport of the defect measures. 
 In section  \ref{weak}, we prove the main theorem \ref{fene}. 
As is classical when proving   global existence of 
weak solutions, the only none trivial part is the proof of  
the weak compactness of a sequence of global solutions satisfying the 
a priori estimates and 
we will only detail  this part of the proof.  In section \ref{app-seq}, 
we present one way of approximating the system. In section \ref{conc} 
we present some concluding remarks and some open problems.

\section{Preliminaries}\label{prem}

\subsection{Hardy type inequalities}
The dissipation  term in the free energy estimate \eqref{free} measures the distance 
between $\psi$ and the equilibrium $\psi_\infty$. We would like 
to use that bound to control the extra stress tensor in $L^2$. 
This  will be done using  the following Hardy \cite{Hardy25}  type inequality. 

\begin{lem} \label{h-lem}

If $k > 1$, then we have 
\begin{equation} \label{Hardy-1}
\int_0^1 \frac{\psi}{x^2}   \leq 
C \int_0^1 x^k \left| \left(\sqrt{\frac{\psi}{ x^k}} \right)  ' \right|^2  + \psi  .  
\end{equation}

For  $k  >   0$,  we have 
\begin{equation} \label{Hardy-inter}
 \left( \int_0^1 \frac{\psi}{x} \right)^2 \leq C  \left(\int_0^1 \psi \right)  \ 
  \left( \int_0^1 x^k | \left(\sqrt{\frac{\psi}{ x^k}} \right)' |^2  + \psi \right)  .  
\end{equation}

 For $ -1 \leq \beta < k \leq 1 $, we have 
\begin{equation} \label{Hardy-inter2}
 \left( \int_0^1 \frac{\psi}{x^{1+\beta}} \right) \leq C 
 \left(\int_0^1  \psi \right)^{1-\beta \over 2}  \ 
  \left( \int_0^1 x^k | \left(\sqrt{\frac{\psi}{ x^k}} \right)' |^2  + \psi \right)^{1+\beta \over 2}  
\quad \quad \hbox{and more generally for all $\gamma \geq 0$}  
\end{equation}

\begin{equation} \label{Hardy-inter-log}
 \left( \int_0^1 \frac{\psi \log^{\gamma} \Big( C+  \frac{\psi}{x^k} 
   \Big)  }{x^{1+\beta}} \right) \leq C  
 \left(\int_0^1  \psi \log^{2 \gamma \over 1 -\beta } \Big( C+   \frac{\psi}{x^k} 
   \Big)    \right)^{1-\beta \over 2}  \ 
  \left( \int_0^1 x^k | \left(\sqrt{\frac{\psi}{ x^k}} \right)' |^2  +
  \psi \right)^{1+\beta \over 2}.   
\quad \quad   
\end{equation}

\end{lem}

\begin{rem} 
Before giving the proof, let us mention that this lemma should be compared to 
the results of section 3.2 of \cite{Masmoudi08cpam}. In particular Proposition 
3.1 was used to control the extra stress tensor. However, the main difference 
is that the results of section 3.2 of \cite{Masmoudi08cpam} are done in an 
$L^2$ frame work since we were dealing with strong solutions there,  however 
the results of lemma \ref{h-lem} are in an $L^1$ frame work since  we only have 
a control on the free energy and its dissipation. 
\end{rem}

Inequality \eqref{Hardy-1}  for  $k> 1$  is just  Hardy inequality. Notice 
that there is no requirement on the  boundary  data since $k > 1$. 
To prove it, we make the change of variable $y = x^{1-k}$ and $ h(y) =  
\sqrt{\frac{\psi(x)}{ x^k}} $. Hence, to prove 
\eqref{Hardy-1}, it is enough to prove that 

\begin{equation} \label{Hardy1-inter-y}
  \int_1^\infty  \frac{h^2}{y^2} dy    \leq C
  \int_1^\infty     h'(y)^2   +  \frac{h^2}{y^{2\alpha}} dy     
\end{equation}
where $\alpha = \frac{k}{k-1} > 1$. To prove \eqref{Hardy1-inter-y}, we integrate by parts 
in 
\begin{equation}
  \int_1^A   \frac{h \, h' }{y} dy    = 
   \int_1^A     \frac{h^2}{ 2 y^{2}}   dy \, + \frac{h(A)^2}{2}   - \frac{h(1)^2}{2}   
\end{equation}
for each $A > 1$. 
The left hand side is bounded by  $  C  (\int_1^A     \frac{h^2}{  y^{2}}   dy)^{1/2} 
   $  $   (\int_1^A     {h'(y)^2}    dy)^{1/2} .$
To bound, $h(1)^2 $ by the right hand side of \eqref{Hardy1-inter-y}, we use that 
$h(y) \leq C \sqrt{y} $ since $ \int_1^\infty     {h'(y)^2}    dy < \infty $ hence, 
$\frac{h^2}{y^\alpha} $ goes to zero when $y$ goes to infinity. This yields that 
\begin{equation}
h^2(1) = -   \int_1^\infty   \left( \frac{h^2 }{y^\alpha} \right)' dy    = 
  -  \int_1^\infty     2  \frac{h}{  y^{\alpha}}  h'  - \alpha \frac{h^2}{y^{\alpha+1}}   dy 
\end{equation}
which is controlled by the right hand side of \eqref{Hardy1-inter-y} 
using Cauchy-Schwarz and the fact that $\alpha > 1$.  Letting $A$ go
to infinity, we get the result.

The proof of \eqref{Hardy-inter}  when $k> 1$ follows by interpolation. 

In the case $0< k \leq 1 $, \eqref{Hardy-1}  only holds if we add a vanishing  boundary 
condition at  $x=0$. However, we can still  prove that  \eqref{Hardy-inter}  holds 
without any   extra condition. Indeed,  making the change of variables  $y = x^{1-k} $ 
(when $k < 1$)
 and denoting $ h(y) =  \sqrt{\frac{\psi(x)}{ x^k}} $, we see that \eqref{Hardy-inter}
is equivalent to

\begin{equation} \label{Hardy-inter-y}
 \left( \int_0^1  y^{\alpha-1}  h^2 dy  \right)^2 \leq C
  \left(\int_0^1   y^{2 \alpha}  h^2\, dy \right)  \ 
  \left( \int_0^1 h'(y)^2   + y^{2 \alpha}  h^2   \right) 
\end{equation}
 where $\alpha = {k \over 1-k}$.  To prove \eqref{Hardy-inter-y}, we 
integrate by parts in the following integral :

\begin{equation}   
\int_0^1  y^{\alpha}  h \  h' dy   = -\frac\alpha2 \int_0^1
 y^{\alpha-1} h^2  + \frac{h^2(1)}{2}.  
\end{equation}
and notice that the left hand side is bounded by 
$  \left(\int_0^1  y^{2\alpha}  h^2 \, \int_0^1 h'(y)^2  \right)^{1/2}  $ using Cauchy-Schwarz
inequality.

Moreover, we have   
\begin{eqnarray} \nonumber   
h(1)^2  = \int_0^1   (y^{2\alpha+ 1 } h^2)' dy &=&  \int_0^1   y^{2\alpha+ 1 } h  \, h' 
 + (2\alpha + 1 ) y^{2\alpha } h^2 \\ 
 \label{h(1)}  &\leq&    C  \left( \int_0^1    h'(y)^{2} +   y^{2\alpha } h^2  
 \,   \int_0^1  y^{2\alpha } h^2  \right)^{1/2}. 
 \end{eqnarray}
Hence,  \eqref{Hardy-inter-y} follows. 

When $k=1$, we make the change of variable $y =-  \log x$ and hence  \eqref{Hardy-inter}
is equivalent to 
\begin{equation} \label{Hardy-inter-y-k=1} 
 \left( \int_0^\infty   e^{-y} h^2 dy  \right)^2   \leq C
  \left(\int_0^\infty   e^{-2y}  h^2\, dy \right)  \ 
  \left( \int_0^\infty  h'(y)^2   +   e^{-2y}  h^2   \right) 
\end{equation}
and the proof of \eqref{Hardy-inter-y-k=1}  can be done in a similar way as that of 
\eqref{Hardy-inter-y}.

To prove \eqref{Hardy-inter2}, we first notice that if $-1 \leq \beta \leq 0 $, then
the inequality can be easily deduced from  \eqref{Hardy-inter} by interpolation. 
 When $\beta > 0$,  \eqref{Hardy-inter2}   is equivalent  (in the case $k<1$)  to 
\begin{equation} \label{Hardy-inter-y2}
 \left( \int_0^1  y^{\alpha_\beta-1}  h^2 dy  \right)^2 \leq C
  \left(\int_0^1  y^{2 \alpha}  h^2\, dy \right)^{1-\beta \over 2}  \ 
  \left( \int_0^1 h'(y)^2   + y^{2 \alpha}  h^2   \right)^{1+ \beta \over 2} 
\end{equation}
where $\alpha_\beta = {k-\beta \over 1 - k}$ and $\alpha = {k \over 1-k}$.
Applying \eqref{Hardy-inter-y} with $\alpha$ replaced by $\alpha_\beta$, we get 

\begin{equation} \label{Hardy-inter-y-alpha-b}
 \left( \int_0^1  y^{\alpha_\beta-1}  h^2 dy  \right) \leq C
  \left(\int_0^1  y^{2 \alpha_\beta}  h^2\, dy \right)^{1/2}  \ 
  \left( \int_0^1 h'(y)^2   + y^{2 \alpha}  h^2   \right)^{1/2} 
\end{equation}
Notice that we kept $\alpha$ in the last term instead of putting
$\alpha_\beta$.
 Indeed, the  last integral comes from 
the estimate of $h^2(1) $ and we can keep  $\alpha = {k \over 1-k}$ in \eqref{h(1)}. 
Now, we can apply \eqref{Hardy-inter-y-alpha-b} replacing $\alpha_\beta -1 $ by 
$2\alpha_\beta$ and we get 
\begin{equation} \label{Hardy-inter-y-alpha-b2}
 \left( \int_0^1  y^{2 \alpha_\beta}  h^2 dy  \right) \leq C
  \left(\int_0^1  y^{2 (2\alpha_\beta+1)}  h^2\, dy \right)^{1/2}  \ 
  \left( \int_0^1 h'(y)^2   + y^{2 \alpha}  h^2   \right)^{1/2} 
\end{equation}
We can iterate this, replacing $\alpha-1$ by $2\alpha_\beta$, $2(2\alpha_\beta+1)$, ...
in \eqref{Hardy-inter-y-alpha-b} till we get an index greater than $2 \alpha  = 2 {k \over 1-k} $. 
Interpolating with the last inequality, yields \eqref{Hardy-inter2}. 

In the  case $k=1$, \eqref{Hardy-inter2} is equivalent to  
\begin{equation} \label{Hardy-inter-y-k=1-beta} 
 \left(   \int_0^\infty   e^{- (1-\beta) y} h^2 dy    \right) \leq C
  \left(   \int_0^\infty   e^{-2y}  h^2\, dy     \right)^{1-\beta \over 2}  \ 
  \left(   \int_0^\infty  h'(y)^2   +   e^{-2y}  h^2   \ dy     \right)^{1+ \beta \over 2} 
\end{equation}
The proof of \eqref{Hardy-inter-y-k=1-beta}    is 
similar  and is left to the reader.

For the proof of \eqref{Hardy-inter-log}, we use that it is equivalent (in the case $k<1$)  to
\begin{equation} \label{Hardy-inter-y-log}
 \left( \int_0^1  y^{\alpha_\beta-1}  h^2 \log^{\gamma} (  h^2)    dy  \right)^2 \leq C
  \left(\int_0^1  y^{2 \alpha}  h^2   \log^{2 \gamma \over 1 -\beta  }  (  h^2)   \, dy \right)^{1-\beta \over 2}   \ 
  \left( \int_0^1 h'(y)^2   + y^{2 \alpha}  h^2   \right)^{1+ \beta \over 2}.  
\end{equation}
Again, one can  prove \eqref{Hardy-inter-y-log} in the case $\beta = 0$ by an integration 
by parts similar to the one used in \eqref{Hardy-inter-y}. 
The case where   $-1 \leq \beta \leq 0 $ can be  deduced by interpolation from  the case $\beta = 0$
and the case $ 0  <  \beta < k     $ can be deduced by a bootstrap argument as the one 
used in the proof of \eqref{Hardy-inter2}.

\subsection{Control of the stress tensor}
We recall that $ \psi_\infty(R)  = \frac{e^{-\U(R)}}{ \int_B e^{-\U(R')} dR'}  
 = (1 - |R|^2)^{k/\beta}  / \int_B (1 - |R'|^2)^{k/\beta} \ dR' $
and since $\beta=1$, $\psi_\infty(R) $  behaves like $(1 -|R|)^k $ when 
$|R|$ goes to 1. In particular we will apply lemma \ref{h-lem} with $x = 1-|R|$.

Using the inequality  \eqref{Hardy-inter} in the radial variable 
 with $x = 1 - |R|$, we get  
 \begin{cor}  \label{tau-est}
There exists a constant $C$ such that 
we have the following bound
\begin{equation}
|\tau(\psi)|^2 \leq (\int_{B} \psi  dR) 
 \int_B  
 \left| \nabla_R  \sqrt{\psi \over \psi_\infty}    \right|^2 \psi_\infty  \,  dR 
\end{equation}
 
\end{cor}

This Corollary can be seen as the $L^1$ version of Proposition 3.1 of \cite{Masmoudi08cpam}. 
It will allow us to control the extra stress tensor by the 
free energy dissipation.  


\subsection{Weighted Sobolev inequality}
In subsection  \eqref{ren-N}, we have to prove the  equi-integrability of $N^n_2$. This 
will require  the control of some higher $L^p$ space of $\sqrt{\psi \over \psi_\infty } $. 
We have the following proposition 

 \begin{prop}  \label{WSI}
There exists $ p > 2  $  and  a constant $C$ such that 
we have the following bound
\begin{equation}\label{WSI-eq}
 \left(  \int_B     \left|  \sqrt{\psi \over \psi_\infty} \right|^p  \psi_\infty  \right)^{1/p} 
\leq 
 \left(   \int_B 
 \left| \nabla_R  \sqrt{\psi \over \psi_\infty}    \right|^2  \psi_\infty    + \psi \,  dR   \right)^{1/2}. 
\end{equation}
\end{prop}

For the proof we first notice that the only difficulty comes from the 
weight and hence we can restrict to the region where $|R| > \frac12$. 
We also use some spherical coordinates, namely $ R = (1-x) \omega $ where 
$\omega \in \S^{D-1}$ and $0< x < \frac12$. The square of the  right hand side of \eqref{WSI-eq}
 can be written as the sum of a radial part and  an angular part :

\begin{equation} \label{rad-p}
  \int_{\S^{D-1}}    \left(   \int_0^{1/2} 
\left[     \left| \partial_x   \sqrt{\psi \over \psi_\infty}    \right|^2  +  |\sqrt{\psi \over \psi_\infty}|^2    \right]  
        \,  x^k \, dx     \right)    d\omega. 
\end{equation}

\begin{equation}\label{ang-p}
 \int_0^{1/2}   \int_{\S^{D-1}}    \left(  \int_{\S^{D-1}}    
\left[     \left| \partial_\omega 
    \sqrt{\psi \over \psi_\infty}    \right|^2  +  |\sqrt{\psi \over \psi_\infty}|^2    \right]  
       d\omega    \right)    \,  x^k \, dx    . 
\end{equation}

We recall the following 1D   weighted  $L^p-L^q$ Hardy inequality (one can also call it 
 weighted Sobolev inequality)

\begin{equation} \label{kmp-th} 
 \left(  \int_0^{1/2}    |F(x)|^q   x^k \, dx  \right)^{1/q}  \leq   
 C    \left(  \int_0^{1/2}    |F'(x)|^2   x^k \, dx  \right)^{1/2}  . 
\end{equation}
This inequality can be easily deduced from Theorem 6 of \cite{KMP07}, taking
$u(x) = v(x) = x^k$ for any $q < \infty$ if $ k \leq 1$ and for 
$q \leq \frac{2(k+1)}{k-1} $ if $k > 1 $. 
Indeed,    Theorem 6 of \cite{KMP07}  stated that    \eqref{kmp-th} holds for 
any $F$, with $F(\frac12) = 0 $ if 
$$  \sup_{0< r < \frac12}   (\int_0^r x^k  dx )^{1/q}   (\int_r^{1 \over 2} (x^k)^{-1}  dx )^{1/2}   
 < \infty.      $$   
Hence, we get a control of 
$  \sqrt{\psi \over \psi_\infty}  $ in the space 
$L^2( \S^{D-1} ; L^q ((0,\frac12),  x^k dx )  ) $ using the 
radial part of the norm  \eqref{rad-p}. 

On the other hand we can use the classical Sobolev inequality in $D-1$ dimension 
to control    $  \sqrt{\psi \over \psi_\infty}  $ in the space  
  $ L^2_x( (0,\frac12) ; L^s (\S^{D-1}) ,   x^k dx    ) $ where 
$s = \frac{ 2(D-1) }{(D-1) -2}  $ if $D> 3$, $s < \infty $ if $D=3$ and 
$s \leq \infty$ if $D=2$. Interpolating between the two spaces 
$L^2_\omega  L^q_x$ and  $L^2_x   L^s_\omega$, we deduce the existence of 
some $p> 2$ such that  \eqref{WSI-eq} holds.

\subsection{Young measures and Chacon limit}

We recall here two important weak convergence objects used in this paper, namely 
the Young measure and the Chacon's biting lemma. Actually, these two notions 
are very related as was observed in  Ball and Murat \cite{BM89}. 

\begin{prop}   (Young measures)    \label{young}
If $f^n$ is a sequence of  functions bounded in $L^1(U;\R^m)$ where $U$ is an open 
set of $\R^N$, then there exists  a family  $(\nu_x)_{x\in U}$ of probability measures 
on $\R^m$ (the Young measures),  
 depending measurably on $x$ and a subsequence  also denoted  $f^n$  such that 
if $g : \R^m \, \to \, \R$ is continuous, if $A \subset U$ is measurable and 
$$  g(f^n)  \rightharpoonup  z(x)  \quad \hbox{weakly in} \, L^1(A; \R),     $$  
then $g(.) \in L^1(\R^m; \nu_x)$  for  a.e.  $x \in A$  and 
$$  z(x) = \int_{\R^m}    g(\lambda) d\nu_x(\lambda) \quad a.e. \quad x \in A.  $$ 
\end{prop}

In the case where $f^n$ is bounded in $L^p(U;\R^m) $ for some $p>1$
(or when $f^n$ is equi-integrable),
we can always take $A = U$ and we have (extracting a subsequence)
$$  g(f^n)  \rightharpoonup  \int_{\R^m}    g(\lambda) d\nu_x(\lambda)   .     $$

\begin{prop}   (Chacon limit)    \label{chac}
If $f^n$ is a sequence of  functions bounded in $L^1(U;\R^m)$ where $U$ is an open 
set of $\R^N$, then there exists  a function $f \in L^1(U;\R^m)  $, a subsequence 
$f^n$ and a non-increasing sequence of measurable sets $E_k$ of $U$ with 
$\lim_{k \to \infty}  \L_N(E_k) = 0 $ (where $\L_N$ is the Lebesgue
measure on $\R^N$) such that for all $k \in \N$, 
 $f^n  \, \rightharpoonup \, f$ weakly in 
$ L^1(U -  E_k;\R^m)  $  as $n$ goes to infinity. $f$ is called the Chacon limit 
of $f^n$. 
\end{prop} 

It is easy to see that if $f^n$ is equi-integrable then the Chacon limit of 
$f^n$ is equal to the weak limit of $f^n$ in the sense of distribution. 

If we consider continuous functions  $g_k :\R^m \, \to \, \R^m $, $k \in \N$ satisfying 
the conditions : 
 
(a) $g_k(\lambda) \to \lambda$ when $k \to \infty$, for each $\lambda \in \R^m$, 

(b)  $| g_k(\lambda) | \leq C ( 1+ |\lambda|  )$, for all $k\in \N$ and $\lambda \in \R^m$, 

(c)  $\lim_{|\lambda| \to \infty}  |\lambda|^{-1} |g_k(\lambda)| = 0
$  for each $k$, 
  
then, under the hypotheses of Proposition \ref{young}, 
 for each fixed $k$, the  sequence of functions $g_k(f^n)$ is equi-integrable and 
hence (extracting a subsequence) converges weakly in  $L^1(U; \R^m)$, to some $f_k$.  
Applying a diagonal process,   
as $k$ goes to infinity, the sequence $f_k$ converges strongly to some 
$f$ in $L^1(U; \R^m)$.  The limit $f$ is the Chacon's limit of the subsequence 
$f^n$ and it is given by 
$$  f(x) = \int_{\R^m}    \lambda  d\nu_x(\lambda) \quad a.e. \quad x \in U .  $$ 
  This gives an other possible  definition of Chacon's limit
    which is equivalent to the one given in Proposition \ref{chac}. For the proof 
of these results we refer to  \cite{BM89}.

\section{A priori estimates}\label{aprio}

\subsection{Free energy}
The second equation of (\ref{micro}) can   be written as 
\begin{equation} \label{psi} 
\partial_t \psi + u. \nabla \psi =   {\rm div}_R \Big[ -  \nabla u  \cdot R \psi\Big] 
      +  div_R  \Big[  \psi_\infty \nabla_R  {\psi \over \psi_\infty} \Big]  .
\end{equation} 
We define $\rho(t,x) = \int_B \psi dR$. Integrating \eqref{psi} in
$R$, we get the transport of $\rho$, namely $\partial_t \rho  +
u. \nabla \rho = 0 .   $

Multiplying \eqref{psi}  by $\log \frac{\psi}{\rho  \psi_\infty}$ and integrating in $R$ and $x$, we get

\begin{equation} \label{psi-est} 
\partial_t  \int_\Omega  \int_B  \psi  \log({\psi \over \rho
  \psi_\infty})  -\psi+ \rho  \psi_\infty 
  =    \int_\Omega    \int_B \nabla u \cdot R \,  \nabla_R \U \psi 
 - 4  \int_\Omega     \int_B  \psi_\infty  
 \left| \nabla_R  \sqrt{\psi \over \psi_\infty}    \right|^2 dR 
\end{equation} 
where we have used that $ \nabla \psi_\infty = - \psi_\infty \nabla \U $.

The first equation of (\ref{micro}) yields the classical energy estimate for 
the Navier-Stokes equation

\begin{equation} \label{u-est} 
\partial_t  \int_\Omega  \frac{|u|^2}{2}
  = -    \int_\Omega    \nabla u : \tau 
 - \nu   \int_\Omega    |\nabla u |^2.   
\end{equation}

Adding \eqref{psi-est}  and \eqref{u-est}  yields the following decay of the 
free-energy

\begin{equation} \label{free-est} 
\partial_t    \int_\Omega   \int_B  [\psi  \log({\psi \over \rho
    \psi_\infty}) - \psi + \rho \psi_\infty    ]
 +   \frac{|u|^2}{2}  
  =     - \nu   \int_\Omega    |\nabla u |^2   
 - 4  \int_\Omega     \int_B  \psi_\infty  
 \left| \nabla_R  \sqrt{\psi \over \psi_\infty}    \right|^2 . 
\end{equation} 

Integrating in time, we get the following uniform bound  for all $t > 0 $

\begin{equation} \label{free} 
 \int_\Omega   \int_B [ \psi  \log({\psi \over \rho
    \psi_\infty}) - \psi + \rho \psi_\infty ]  
 +   \frac{|u|^2}{2}  \ (t) + 
 \int_0^t    \nu   \int_\Omega    |\nabla u |^2   
  +  4  \int_\Omega     \int_B  \psi_\infty  
 \left| \nabla_R  \sqrt{\psi \over \psi_\infty}    \right|^2   =  C_0. 
\end{equation}

To simplify the notations  in the rest of this section, we will 
assume that $\rho_0(x)  = 1$. The proof in the general case is
identical and we will indicate the changes to be made at the end. 
The general idea is the following:  When proving a priori estimates,  
one has just to replace $\psi_\infty$ by $\rho(t,x) \psi_\infty$
and take advantage of the fact that $\rho$ is just transported by the 
flow.  When proving weak compactness,  one can    use that  
$\rho^n$ converges strongly to $\rho$ in all $L^p((0,T)\times \Omega )
$ spaces and use  $\rho^n(t,x) \psi_\infty$. Due to the local character of 
the proof of weak compactness,  a  simpler way  is to just  use 
 $\psi_\infty $ and so the calculations  given in section \ref{weak} 
hold even when $\rho_0$ is not constant.

\subsection{$\log^2$ estimate}
The free energy only gives an  $L\log L  (\psi_\infty dR)  $ bound on $\frac{\psi}{\psi_\infty}$. 
For some technical reasons, we will need to control a slightly higher 
growth of $\psi$ in the $R$ variable.

We introduce $\tilde \psi = \psi + a  \psi_\infty$ for some $a  > 1$. 
This is done to insure that $\log \frac{\tilde \psi}{\psi_\infty} $ does not 
take negative values. It will also add a new term in the equation which will not present any 
extra difficulties.   Hence, $\tilde \psi$ solves 
\begin{equation} \label{psi-a} 
\partial_t \tilde \psi + u. \nabla \tilde \psi =   {\rm div}_R \Big[ -  \nabla u  \cdot R \tilde \psi\Big] 
      +  div_R  \Big[  \psi_\infty \nabla_R  {\tilde \psi \over \psi_\infty} \Big] 
  -  a  \nabla  u \cdot R \psi_\infty \nabla_R \U     .
\end{equation}

We first derive this extra bound in the case the domain 
$\Omega$ is bounded and then discuss the modification of the argument 
in the whole space  case. 
\subsubsection{Case of a  bounded domain} 
Multiplying \eqref{psi-a} 
 by $\log^2 \frac{\tilde \psi}{\psi_\infty}$ and integrating by parts in $R$, we get 
 
  \begin{eqnarray}\nonumber 
& \displaystyle   (\partial_t + u. \nabla_x)     \int_B  \tilde \psi   
[\log^2({\tilde \psi \over \psi_\infty})  -2 \log({\tilde \psi \over \psi_\infty})  + 2 ] 
 =    -2 a k  \nabla_i u_j \int_B \frac{R_i R_j}{1-|R|^2} \psi_\infty  
 \log^2 \frac{\tilde \psi}{\psi_\infty}
  \qquad\\
  \label{log^2}  
& \qquad \qquad \qquad \qquad \displaystyle           +   \int_B \nabla u \cdot R \tilde \psi  \,  2\log({\tilde \psi \over \psi_\infty}) 
   {\psi_\infty \over \tilde \psi}    \nabla_R  {\tilde \psi \over \psi_\infty}  \, -8  
       \int_B  \psi_\infty  
 \left| \nabla_R  \sqrt{\tilde \psi \over \psi_\infty}    \right|^2 \log({\tilde \psi \over \psi_\infty})     
\end{eqnarray} 
   
The second  term on the right hand side of \eqref{log^2} can be rewritten

\[
2     \int_B \nabla u \cdot R \psi_\infty  \,  
\nabla_R  \left( {\tilde \psi \over \psi_\infty}\log({\tilde \psi \over \psi_\infty})  -  {\tilde \psi \over \psi_\infty}    \right)    = 2      \int_B \nabla u \cdot R \, \nabla_R \U 
 \tilde \psi \left(  \log({\tilde \psi \over \psi_\infty}) -1 \right)
\]

 Taking the square root of \eqref{log^2}, we get

 \begin{eqnarray}\nonumber 
& \displaystyle   (\partial_t + u. \nabla_x)    \left(  \int_B  \tilde \psi   
[\log^2({\tilde \psi \over \psi_\infty})  -2 \log({\tilde \psi \over \psi_\infty})  + 2 ] \right)^{1/2} 
 =    \frac{ -  a k  \nabla_i u_j \int_B \frac{R_i R_j}{1-|R|^2} \psi_\infty  
 \log^2 \frac{\tilde \psi}{\psi_\infty} } {  \left( \int_B  \tilde \psi   
[\log^2({\tilde \psi \over \psi_\infty})  -2 \log({\tilde \psi \over \psi_\infty})  + 2 ] \right)^{1/2}  } \qquad \\
 \label{sqrt-log^2}  
& \qquad \qquad \qquad \qquad +  \displaystyle      
  \frac{  \int_B \nabla u \cdot R \tilde \psi  \,  2\log({\tilde \psi \over \psi_\infty}) 
   {\psi_\infty \over \tilde \psi}    \nabla_R  {\tilde \psi \over \psi_\infty}  }{\left(  \int_B  \tilde \psi   
[\log^2({\tilde \psi \over \psi_\infty})  -2 \log({\tilde \psi \over \psi_\infty})  + 2 ] \right)^{1/2}  } 
-4 
   \frac{     \int_B  \psi_\infty  
 \left| \nabla_R  \sqrt{\tilde \psi \over \psi_\infty}    \right|^2 \log({\tilde \psi \over \psi_\infty}) }
  {\left(  \int_B  \tilde \psi   
[\log^2({\tilde \psi \over \psi_\infty})  -2 \log({\tilde \psi \over \psi_\infty})  + 2 ] \right)^{1/2}  }  \\ 
& \nonumber  \qquad \qquad \qquad \qquad  = I_1 + I_2 + I_3.       
\end{eqnarray} 

Let us introduce the notation 

\begin{equation}  \label{N_2-def} 
 N_2 =   \left(  \int_B  \tilde \psi   
[\log^2({\tilde \psi \over \psi_\infty})  -2 \log({\tilde \psi \over \psi_\infty})  + 2 ] \right)^{1/2}.   
\end{equation}

To bound 
$I_1$ we use that,  
$ \psi_\infty  
 \log^2 \frac{\tilde \psi}{\psi_\infty}  \leq C \tilde \psi  $. Hence, the numerator 
of $I_1$ is bounded by $ C  |\nabla u | \int \frac{\tilde \psi}{1-|R|^2} dR   $ 
which is clearly in $L^1((0,T)\times \Omega \times B ) $. Indeed, by
using \eqref{Hardy-inter} and  Corollary \eqref{tau-est}, 
we see that 

\begin{equation}  
(\int \frac{ \psi}{1-|R|^2} dR) \leq C   \left(  \int_B  \psi_\infty  
 \left| \nabla_R  \sqrt{\psi \over \psi_\infty}    \right|^2 dR \right)^{1/2} .  
\end{equation}

To bound the second  term on the right hand side of \eqref{sqrt-log^2}, we use that 
the numerator can be bounded by 

 \begin{eqnarray}\nonumber  
& &  \!\!\!\!\!\!\!\!\!\! \!\!\!\!\!\!   \displaystyle  \left|    \int_B \nabla u \cdot R \tilde \psi  \,  \log({\tilde \psi \over \psi_\infty}) 
   {\psi_\infty \over \tilde \psi}    \nabla_R  {\tilde \psi \over \psi_\infty}     \right|  \leq  \qquad \qquad  \qquad \qquad \\
  \label{Num-bound}  
& 
 \leq & C   |\nabla u| \left( \int_B   \psi_\infty  
  | \log({\tilde \psi \over \psi_\infty}) |  
     \left|  \nabla_R  \sqrt{\tilde \psi \over \psi_\infty} \right|^2     \right)^{1/2}
\left( \int_B   \psi_\infty  
  | \log({\tilde \psi \over \psi_\infty}) |  
      {\tilde \psi \over \psi_\infty}      \right)^{1/2} \\
&  
 \leq  & C      |\nabla u|^2  \left( \int_B   \tilde \psi  
  | \log({\tilde \psi \over \psi_\infty}) |  
           \right)  +   \left( \int_B   \psi_\infty  
  | \log({\tilde \psi \over \psi_\infty}) |  
     \left|  \nabla_R  \sqrt{\tilde \psi \over \psi_\infty} \right|^2     \right) \\  
& 
 \leq& C       |\nabla u|^2   (1+a)^{1/2}   \left( \int_B   \tilde \psi  
   \log^2({\tilde \psi \over \psi_\infty})   
           \right)^{1/2}  +   \left( \int_B   \psi_\infty  
  | \log({\tilde \psi \over \psi_\infty}) |  
     \left|  \nabla_R  \sqrt{\tilde \psi \over \psi_\infty} \right|^2     \right)  
\end{eqnarray} 
  Dividing \eqref{Num-bound}  by $N_2 $, we deduce that 

\begin{equation}  
 I_2 \leq  C  |\nabla u|^2 - \frac14 I_3.   
\end{equation}


Integrating \eqref{sqrt-log^2} in time and space and using the fact that 
$I_3$ has a sign, we deduce the following a priori bound

  \begin{equation} \label{log2-b}
  \int_\Omega  \left(  \int_B  \tilde \psi   
[\log^2({\tilde \psi \over \psi_\infty})  -2 \log({\tilde \psi \over \psi_\infty})  + 2 ] \right)^{1/2}  (t) + \int_0^T  \int_\Omega
   \frac{     \int_B  \psi_\infty  
 \left| \nabla_R  \sqrt{\tilde \psi \over \psi_\infty}    \right|^2 \log({\tilde \psi \over \psi_\infty}) }
  {\left(  \int_B  \tilde \psi   
[\log^2({\tilde \psi \over \psi_\infty})  -2 \log({\tilde \psi \over \psi_\infty})  + 2 ] \right)^{1/2}  }  \leq C_T
\end{equation} 
for $0\leq t \leq T$, if the initial condition satisfies 
$   \int_\Omega  \left(  \int_B  \tilde \psi_0   
[\log^2({\tilde \psi_0 \over \psi_\infty})  -2 \log({\tilde \psi_0 \over \psi_\infty})  + 2 ] \right)^{1/2}      \leq C_0. $
 Hence, we see that \eqref{sqrt-log^2} can be written as 
 \begin{equation} \label{NF}
(\partial_t + u.\nabla) N_2 = F_2  
\end{equation} 
where $F_2$ is in $L^1((0,T) \times \Omega)$. 

It turns out that passing to the limit in the bound  \eqref{log2-b} is 
not clear. Actually, one can find sequences of functions $\tilde \psi^n$ such that 
\eqref{log2-b} holds and the weak limit does not satisfy 
\eqref{log2-b}. This is the reason, we prefer to write the second bound as

 \begin{equation} \label{log2-b2}
  \int_\Omega  \left(  \int_B   g^2 \log (g^2) \psi_\infty  dR      \right)^{1/2} dx +  \int_0^T  \int_\Omega
   \int_B   \frac{     \psi_\infty  
 \left| \nabla_R  g   \right|^2    \ dR}
  {\left(  \int_B  \tilde \psi   
[\log^2({\tilde \psi \over \psi_\infty})  -2 \log({\tilde \psi \over \psi_\infty})  + 2 ] \right)^{1/2}  }  \leq C_T
\end{equation} 
where $g$ is given by $g =  \sqrt{{\tilde \psi \over \psi_\infty}} 
  \log^{1/2}({\tilde \psi \over \psi_\infty}) $. 
 
\subsubsection{Case of unbounded domain} 
 
In the case $\Omega = \R^D $, we first take $c_1$ and $c_2$ 
the two constants  such that the function 
$ \phi(x)=  x[\log^2 x -2 \log x + c_1  ] +c_2  $ satisfies 
 the fact that $\phi(1+a) = \phi'(1+a ) = 0   $. This is 
achieved by taking $c_1 = 2 - \log^2(1+a)$ and 
$c_2 = 2 (1+ a ) [\log(1+a) - 1   ] $. Notice also 
that the function $\phi (x)  $ is nonnegative for $x \geq a $ 
 since $a$ is taken big enough. 
It is clear that the extra  bound \eqref{ini-log2}
implies that 
\begin{equation}
  \int_\Omega     \frac{  \int_B \phi( {\tilde \psi_0  \over \psi_\infty}
    )  dR    }{ 1+  \left[  
  \int_B \phi( {\tilde \psi_0  \over \psi_\infty} )  dR  \right]^{1/2}  }     dx \leq C_0   
\end{equation} 
 and hence, we can perform the same calculations as \eqref{log^2}  and
 \eqref{sqrt-log^2} with  
$  \int_B  \tilde \psi   
[\log^2({\tilde \psi \over \psi_\infty})  -2 \log({\tilde \psi \over
    \psi_\infty})  + 2 ]  dR    $ replaced by 
$  \int_B  \phi( {\tilde \psi \over \psi_\infty}
    )  dR    $  and with the function $s \to \sqrt{s}$ used to go from 
\eqref{log^2}  to 
 \eqref{sqrt-log^2} replaced by $ s \to \frac{s}{1 + \sqrt{s}}  $ 
which behaves like  $\phi_1(s) = $   min$(\sqrt{s},s)$. The rest of the proof is
identical.  

\subsubsection{Case $\rho$ is not constant} 
In the case  $\rho$ is not constant and we are in a bounded domain,
 we have to modify
\eqref{log^2} slightly 
and multiply by  $\log^2\frac{\tilde \psi}{\rho \psi_\infty}  $. 
In the case we are also in an unbounded domain, we have to replace 
$  \int_B  \tilde \psi   
[\log^2({\tilde \psi \over \psi_\infty})  -2 \log({\tilde \psi \over
    \psi_\infty})  + 2 ]  dR    $  by 
$  \int_B  \phi(  {  (1+a )  \tilde \psi \over   (\rho + a )   \psi_\infty}
    )  dR    $. The extra factor $\frac{1+a}{\rho+a}$ is used to 
insure that when  $\tilde \rho$ is at microscopic equilibrium, namely  
$\tilde \psi = (\rho + a ) \psi_\infty$, the integrand reduces to 
$\phi(1+a)$.  The rest of the proof is identical and yields at the end 
the following bound instead of \eqref{log2-b}

 \begin{equation} \label{log2-b-mod}
  \int_\Omega   \phi_1 \left(  \int_B   \phi({  (1+a )  \tilde \psi
    \over   (\rho + a )   \psi_\infty}   )   \right)  (t) + \int_0^T  \int_\Omega
   \frac{     \int_B  \psi_\infty  
 \left| \nabla_R  \sqrt{\tilde \psi \over \psi_\infty}    \right|^2
 \log({ (1+a )  \tilde \psi \over  (\rho +a ) \psi_\infty}) }
  { 1 +  \left(  \int_B   \phi({  (1+a )  \tilde \psi
    \over   (\rho + a )   \psi_\infty}   )   \right)^{1/2}        }
  \leq C_T. 
\end{equation} 
One can then deduce from \eqref{log2-b-mod} that \eqref{log2-b} and \eqref{log2-b2} hold 
with  the integration set $\Omega$ replaced by any compact $K$ of
$\R^D$.

\section{Weak compactness} \label{weak}
As  it is classical when proving global existence of weak solutions,
 it is enough  to prove the weak 
compactness of a sequence of weak solutions satisfying the a priori estimates of the 
previous section. In the next section, we present one way of approximating 
the system.  We consider  $(u^n, \psi^n)$ a sequence of weak solutions to  
\eqref{micro} satisfying, uniformly in $n$,  the free energy bound  \eqref{free}   and the $\log^2$ bound 
\eqref{log2-b}  with an initial data  $(u^n_0, \psi^n_0)$  such that 
 $(u^n_0, \psi^n_0)$  converge strongly to  $(u_0, \psi_0)$ in $L^2(\Omega) \times 
 L^1_{loc}(\Omega; L^1( B))$ 
 and $ \psi^n_0 \log \frac{\rho^n_0 \psi^n_0 }{\psi_\infty } - \psi^n_0 + \rho^n_0 \psi_\infty $
  converges strongly to $\psi_0  \log \frac{\psi_0}{\rho_0 \psi_\infty} - \psi_0 + \rho_0 \psi_\infty $
in $  L^1(\Omega \times B) $. We also assume that $(u^n, \psi^n) $ has
some 
extra regularity with bounds that depend on $n$ such that we can
perform  all the following calculations. 

We extract a subsequence such that $u^n$ converges weakly to $u$
in $L^p((0,T); L^2(\Omega)) \cap  L^2((0,T); H^1_0(\Omega))$ and $\psi^n$ 
converges weakly to $\psi$ in   $L^p((0,T) ; L^1_{loc} (  \Omega \times B  ))$
for each $p< \infty$. We would like to prove that $(u,\psi)$ is still a solution of 
\eqref{micro}. The main difficulty is to pass to the limit in the nonlinear term 
$\nabla u^n   R \psi^n $ appearing in  the second equation of \eqref{micro}. 
 
We  introduce $g^n = \sqrt{\frac{\tilde \psi^n}{\psi_\infty}}
\log^{1/2}({\tilde \psi^n  \over \psi_\infty})  $    and $f^n = \sqrt{\frac{\tilde \psi^n}{\psi_\infty}}$
where $ \tilde \psi^n =   {\psi^n + a \psi_\infty}$ and 
$a > 1 $ is any constant.  We also assume, extracting a subsequence if necessary,  that
$g^n$ and $f^n$  converge weakly to some $g$  and $f$ 
in  $L^p((0,T) ; L^2_{loc}(  \Omega \times B, \, dx \psi_\infty dR  ))$
for each $p< \infty$.  To prove that $(u,\psi)$ is  a solution of  
\eqref{micro}, it will be enough to prove that 
$ (g^n)^2 = \frac{\tilde \psi^n}{\psi_\infty} 
\log({\tilde \psi^n \over \psi_\infty})   $ converges weakly to 
$g^2= \frac{\tilde \psi}{\psi_\infty} 
\log({\tilde \psi \over \psi_\infty})  $, namely that 
$g^n$ converges strongly to $g$ in $L^2((0,T) ; L^2(  \Omega \times B, \, dx \psi_\infty dR 
))$.

First, it is clear that $u,\tilde \psi$ and $ g$ satisfy  the same a priori estimates that
the sequence   $u^n,\tilde \psi^n$ and $ g^n$ satisfy since all those 
functionals have good convexity properties.  In particular
it is clear that $u,\psi$ satisfy \eqref{free}. We just point out that
to pass to the limit in the last term on the left hand side of
\eqref{free}, 
we can use the fact that the function   $ \phi_2  (x,y) =  \frac{x^2}y $ is
convex.  To pass to the limit in   \eqref{log2-b2}, we also use 
the fact that  $ \phi_2  (x,y) $ is convex. Hence, we deduce that  
 \begin{equation} \label{log2-b-lim}
\Sup_{0\leq t \leq T}  \int_\Omega  \left(   \Big(  \int_B   g^2 \log (g^2) \psi_\infty      dR      \Big)^{1/2} + 
   \o{N^n_2} \right)  dx  (t)  + \int_0^T  \int_\Omega
   \int_B   \frac{     \psi_\infty  
 \left| \nabla_R  g   \right|^2 }
  {   \o{N^n_2}  }  \leq C_T 
\end{equation} 
where $\o{N^n_2} $ is the weak limit of 
$      
\left(  \int_B  \tilde \psi^n   
[\log^2({\tilde \psi^n \over \psi_\infty})  -2 \log({\tilde \psi^n  \over
    \psi_\infty})  + 2 ] \right)^{1/2}  $.  

Dividing  \eqref{psi-a} by $\psi_\infty$ we get 
\begin{equation} \label{psi-a-d} 
\partial_t   \frac{\tilde \psi^n}{\psi_\infty} + u^n. \nabla \frac{\tilde \psi^n}{\psi_\infty}    =
   {\rm div}_R \Big[ -  \nabla u^n  \cdot R     \frac{\tilde \psi^n}{\psi_\infty}  \Big] 
      + \nabla \U. \nabla u^n R  \frac{\tilde \psi^n}{\psi_\infty}  + \Delta_R \frac{\tilde \psi^n}{\psi_\infty} 
        - \nabla \U. \nabla_R \frac{\tilde \psi^n}{\psi_\infty}  
  -  a  \nabla  u^n \cdot R  . \nabla_R \U     .
\end{equation}

From \eqref{psi-a-d}, we deduce that for any   smooth function 
   $\Th$ from $(0,\infty)$ to $\R$,   we have

\begin{equation} \label{g^n} 
\begin{split}  
\partial_t \Th(\pp)  + u^n. \nabla \Th(\pp) &  =  
 -  \nabla u^n R  \cdot  \nabla_R \Th(\pp)     
 + \nabla_R \U . \nabla u R \pp \Th'(\pp) \\
& \quad +  \Delta \Th(\pp)  - \nabla_R \U \cdot \nabla_R \Th(\pp) 
- \Th''(\pp ) |\nabla_R \pp|^2 \\ 
 & \quad   - 2   ak   \nabla_i  u^n_j   \frac{R_i R_j}{1-|R|^2}  \Th'(\pp)    .
\end{split}
\end{equation}

We take  $\Th(t) = t^{1/2} \log^{1/2} (t) $ and  recall  that 
$g^n = \Th(\pp)$. We introduce the following defect measures 
$\gamma_{ij}, \gamma'_{ij} $ and $\beta_{ij}$ such that 
\begin{equation} \label{defect1}
 \begin{split}  
\nabla  u^n g^n  \to \nabla u g + \gamma,   \quad  \quad  
\nabla  u^n   \pp \Th'(\pp) 
   \to \nabla u  \overline{ \pp \Th'(\pp)     } + \gamma' 
\\ 
\nabla u^n \tilde \psi^n  \to \nabla u \tilde \psi + \beta 
 \quad  \quad  \quad  \quad  \quad  \quad  \quad  \quad  \quad  \quad  \quad  \quad  
 \end{split}  
\end{equation}
where  $\gamma, \gamma'  \in L^2((0,T) \times \Omega \times B)$  and  
$\beta  \in L^2((0,T) \times \Omega ; L^1 (B))$  are matrix valued.

On one hand, passing to the limit in \eqref{g^n} with 
 $\Th(t) = t^{1/2} \log^{1/2} (t)  $, we get

\begin{equation} \label{f-eq}
\begin{split}  
\partial_t g   + u. \nabla g  &  = 
   {\rm div}_R \Big[ -  \nabla_i  u_j   R_j  g   -\gamma_{ij} R_j
   \Big]  \\ 
& \quad + \nabla_R \U . \nabla u R \overline{\pp \Th'(\pp) } 
   +  \nabla_R \U  R  : \gamma' \\
& \quad + \frac1{\psi_\infty}  {\rm div}_R \Big[  \psi_\infty \nabla_R g  \Big]  
 +  \overline{ \frac{|\nabla_R f^n|^2 (\log^{1/2} + \log^{-3/2})(\pp)}{f^n}}   \\ 
 & \quad   -   ak \,   \overline {\Th'(\pp) \nabla_i u^n_j   }  
   \frac{2 R_i R_j}{1-|R|^2}     
\end{split}
\end{equation}   
 where, here and below, $\overline{F_n}$ denotes the weak limit of $F_n$ and where 
 we have used that

\begin{equation} \label{Th-deriv}
\left\{   \begin{split}  
\Th' (s) & =  \frac12 s^{-1/2}  ( \log^{1/2} (s)  + \log^{-1/2} (s)   )  \\
\Th'' (s) &   =  -\frac14 s^{-3/2}  ( \log^{1/2} (s)  + \log^{-3/2} (s)   ).   
\end{split}  \right.
\end{equation}

Multiplying by $g$, we get 
\begin{equation} \label{g2-eq}
\begin{split}  
\partial_t g^2   + u. \nabla g^2  &  =   {\rm div}_R \Big[ -  \nabla u
    R g^2\Big]
+ \nabla uR \cdot \nabla_R \U   \overline{ \Big( \pp \Th'(\pp) \Big)  }
   2g  \\ 
& \quad -   {\rm div}_R  (\gamma_{ij} R_j) 2g  +  \nabla_R \U  R  : \gamma'2g \\
 & \quad + \frac1{\psi_\infty}  {\rm div}_R \Big[  \psi_\infty
  \nabla_R g^2  \Big]      - 2 |\nabla_R g|^2  \\ 
& \quad  +  \overline{ \frac{|\nabla_R f^n|^2 (\log^{1/2} + \log^{-3/2})(\pp) 
     |^2}{f^n}} 2g   
  -   ak \,   \overline {\Th'(\pp) \nabla_i u_j^n   }  
   \frac{4 R_i R_j}{1-|R|^2}  g    
\end{split}
\end{equation}   

Multiplying \eqref{g2-eq} by ${\psi_\infty} $
 and integrating in $R$ yields

\begin{equation} \label{g2log-eq}
\begin{split}    
&  (\partial_t  + u. \nabla )  \int_B \psi_\infty g^2 
        =  -  \nabla u : \tau \left(\psi_\infty \Big(   g^2 - 2g \overline{\pp
     \Th'(\pp) } \Big)   \right)  \\
& \quad - \int_B  
   {\rm div}_R  (\gamma_{ij} R_j) 2g \psi_\infty  +  \nabla_R \U  R  :
   \gamma'2g \psi_\infty \\ 
& \quad  +  \int_B   \psi_\infty   \overline{ \frac{|\nabla_R f^n|^2 (\log^{1/2} + \log^{-3/2})(\pp) 
     |^2}{f^n}} 2g          -2 \psi_\infty |\nabla_R g|^2    \\
& \quad   -  \int_B  \psi_\infty  ak \,   \overline {\Th'(\pp)
     \nabla_i u_j^n    }  
   \frac{4 R_i R_j}{1-|R|^2}  g 
\end{split}
\end{equation}

where  we recall that  $\tau_{ij} (\psi) = 2k \int_B \psi 
  \frac{R_i R_j}{1-|R|^2} dR $. 
Here,  there is a small problem of definition:  The  terms on the
second line of the  right hand side 
are  not well defined in the sense of distribution and we need some 
further analysis to make sense of them.  Also, the transport term 
is not well defined even if we write it as div$(u \int_B \psi_\infty
g^2)$.  Actually, as we
will see later,   we 
will not use \eqref{g2log-eq} but a renormalized form of it. 
 Indeed, we will construct in the next subsection 
a renormalizing factor $N$ that satisfies $(\partial_t  + u. \nabla )
\frac1N =0 $  and we will make sense of \eqref{g2log-eq} after
dividing 
each term by $N^4$.

On the other hand, passing to the limit in the equation satisfied by $\tilde \psi_n$, we get 

\begin{equation} \label{psi-a-lim} 
\partial_t \tilde \psi + u. \nabla \tilde \psi =   
{\rm div}_R \Big[ -  \nabla u  \cdot R \tilde \psi  - \beta_{ij} R_j  \Big] 
      +  div_R  \Big[  \psi_\infty \nabla_R  {\tilde \psi \over \psi_\infty} \Big] 
  -   2 a  k  \nabla  u  :   \psi_\infty   \frac{R_iR_j}{1-|R|^2}    .
\end{equation}   

Besides, $\tilde \psi^n \log( \frac{\tilde \psi^n}{\psi_\infty} )  $ satisfies

\begin{equation} \label{psi-n} \begin{split}    
& (\partial_t + u^n. \nabla)    \left[  \int_B    \tilde 
 \psi^n  \log( \frac{\tilde \psi^n}{\psi_\infty} )     \right]    = 
   \nabla u^n : \tau (\tilde \psi^n) \\
 & \quad \quad \quad  -  4 
 \int_B \psi_\infty \left|    \nabla_R  \sqrt{\frac{\tilde \psi^n}{{\psi_\infty}}}     \right|^2 
  -  2ak   \int_B   \nabla  u^n    \frac{R_iR_j}{1-|R|^2} \psi_\infty 
 \log( \frac{\tilde \psi^n}{\psi_\infty} )      .
\end{split}    
\end{equation}   

We would like  to pass to  the limit weakly in \eqref{psi-n} and deduce that

\begin{equation} \label{psi-n-limit-no} 
\begin{split} 
(\partial_t + u. \nabla)    \int_B    
\overline{\tilde   \psi^n  \log( \frac{\tilde \psi^n}{\psi_\infty} )     }    &  = 
   \nabla u : \tau (\tilde \psi)  + \int_B \beta_{ij} \frac{R_i R_j}{1-|R^2|}  \\ 
&  -  \int_B \psi_\infty   \overline 
{ \left|    \nabla_R  \sqrt{\frac{\tilde \psi^n}{{\psi_\infty}}}     \right|^2 } 
  -  2ak \int_B    \overline{\log( \frac{\tilde \psi^n}{\psi_\infty} )    \nabla  u^n }    \frac{R_iR_j}{1-|R|^2} \psi_\infty  .
\end{split}
\end{equation}   
However,  we can not use \eqref{defect1} to pass 
to the limit in $ \nabla u^n : \tau (\tilde \psi^n)  = \int_B  
  \nabla u^n \frac{R_i R_j }{1 - |R|^2} \tilde \psi^n    $ and deduce that 

\begin{equation} \label{weak-u-t}
 \nabla u^n : \tau (\tilde \psi^n)    \rightharpoonup  \nabla u : \tau (\tilde \psi)  + \int_B \beta_{ij} \frac{R_i R_j}{1-|R^2|} 
\end{equation}   
since  $ \nabla u^n \frac{R_i R_j }{1 - |R|^2} \tilde \psi^n    $ is only bounded 
in $L^1(dt\, dx\, dR  )$. Besides, we can not pass to the limit in the
 transport term even if we write it in divergence form.

To overcome these difficulties, we will divide  \eqref{psi-n} by  $1+\delta N_2^n $ where  
$N_2^n $ 
 solves \eqref{NF}  before 
passing to the limit. Then, we will send $\delta$ to zero. 
 
To be able to deal with the limit $\delta$ to zero,  we need to renormalize \eqref{psi-n} too. 
We denote $N_1^n =  \int_B    \tilde 
 \psi^n  \log( \frac{\tilde \psi^n}{\psi_\infty} )   $, 
$N_2^n = \left(  \int_B  \tilde \psi^n   
[\log^2({\tilde \psi^n  \over \psi_\infty})  -2 \log({\tilde \psi^n \over \psi_\infty})  + 2 ] \right)^{1/2}    $     and 
$ \th_\ka  (s) = \frac{s}{1+\ka s} $.
We first  multiply  \eqref{psi-n} by $\th_\ka'(N_1^n)  $ and get an equation for 
$\th_\ka(N_1^n)  $. Dividing  
 the resulting equation  by  $1+ \delta N_2^n $,  using \eqref{NF} and taking the weak limit when 
$n$ goes to infinity (extracting a subsequence if necessary), we get  for $\ka, \delta > 0$
 

\begin{align}
\nonumber (\partial_t + u . \nabla)  \overline  { \frac{\th_\ka(N_1^n)}{{1 + \delta N_2^n }  }     } 
    = & 
  \overline {  \nabla u^n : \frac{\tau (\tilde \psi^n) }{(1 + \delta N_2^n )(1 + \ka N_1^n )^2 } }  - 
  \overline {  \frac1{(1 + \delta N_2^n ) (1 + \ka N_1^n)^2 } \int_B \psi_\infty \left|    \nabla_R  \sqrt{\frac{\tilde \psi^n}{{\psi_\infty}}}     \right|^2  }  \\ 
&   \label{psi-n-delta-De}  -  \overline {  \frac{2ak}{(1 + \delta N_2^n )(1 + \ka N_1^n)^2 }  
\int_B   \nabla  u^n     \frac{R_iR_j}{1-|R|^2} \psi_\infty 
 \log( \frac{\tilde \psi^n}{\psi_\infty} ) }  \\
&  \nonumber   -  
 \overline { \frac{\delta F^n}{ ( 1 + \delta N_2^n   )^2} {\th_\ka(N_1^n)}  }  .
\end{align}

Now, we can send $\delta $ to zero. Notice that due to the fact that ${\th_\ka(N_1^n)}  $ is 
bounded and that $ F^n  $ is bounded in $L^1 $, we deduce that the last term goes 
to zero when $\delta  $ goes to zero.  Then, we send $\ka$ to zero and recover at 
the limit 

\begin{align}
\nonumber (\partial_t + u . \nabla) \th  
    = & 
  \overline {  \nabla u^n : {\tau (\tilde \psi^n) } }^{\de,\ka}  - 
  \overline {   \int_B \psi_\infty \left|    \nabla_R  \sqrt{\frac{\tilde \psi^n}{{\psi_\infty}}}     \right|^2  }^{\de,\ka}  \\ 
&   \label{psi-n-delta}  -  \overline {  {2ak}  
\int_B   \nabla  u^n     \frac{R_iR_j}{1-|R|^2} \psi_\infty 
 \log( \frac{\tilde \psi^n}{\psi_\infty} ) }^{\de,\ka} . 
\end{align}    
where $\th = \lim_{\ka \to 0} \lim_{\de \to 0 }  \overline  { \frac{\th_\ka(N_1^n)}{{1 + \delta N_2^n }  }}  
  = \lim_{\ka \to 0}  \overline  { {\th_\ka(N_1^n)} }   $  is the Chacon limit of $N_1^n$ 
and  
$$  \overline { F_n}^{ \de,\ka}   = \lim_{\ka \to 0} \lim_{\de \to 0 } 
  \overline {\frac{F_n}{(1 + \delta N_2^n )(1 + \ka N_1^n)^2  } }    $$
for any sequence $F_n$ bounded in $L^1$. We will prove in the next subsection 
that $N^n_1$ is equiintegrable and hence that $\theta$ is also the weak limit of 
$N^n_1$. 
We also notice that if $ F_n $ is equi-integrable then $  \overline { F_n}^{ \de,\ka} = \overline { F_n}$.

To be really precise  the term $u.\nabla \th = \div(\th u )$ on the 
 left hand side of \eqref{psi-n-delta}  is not well defined  since $\th$ is only in 
$L^\infty_t L^1_x$ and $u$ is in $L^2_t  \dot H^1_x$.
 To give sense to \eqref{psi-n-delta}  we need to use a renormalizing 
factor. Actually, we made   the same remark   for the transport term
in   the equation 
\eqref{g2log-eq}. Recall that at the end, we would like to take the difference 
between \eqref{psi-n-delta} and \eqref{g2log-eq} after dividing it by 
the factor $N^4$ that we are going to define now.  The fact of dividing 
by $N$ will insure that $\frac{\theta}N$ is bounded and that all the
terms 
will make sense.   So the point is to divide \eqref{psi-n-delta-De}  by $N^4$ and 
then send $\de$ and $\ka$ to zero.


\subsection{The renormalizing factor $N$} \label{ren-N}
 

We recall that from, \eqref{NF}, $\beta_M(N_2^n)$ solves 
  \begin{equation} \label{NF-renor}
(\partial_t + u^n.\nabla) \beta_M(N_2^n)  = \beta_M'(N_2^n)   F^n .  
\end{equation}  
where $\be_M(s) = \th_{1/M}(s) = \frac{s\, M }{M+s} $. 
Passing to the limit in \eqref{NF-renor}, we get 
 \begin{equation} \label{NF-renor-li}
(\partial_t + u.\nabla)  \overline{\beta_M(N_2^n)}  = \overline{ \beta_M'(N_2^n)   F^n}  .  
\end{equation} 
This equation does not  seem to be very useful since the right hand side is a measure. 
To overcome this problem, we   first  introduce the 
unique a.e flow $X^n$  in the sense of DiPerna and Lions \cite{DL89im}
of $u^n$,  solution of 
\begin{equation} 
\partial_t X^n(t,x) = u^n (t,  X(t,x)  )  \quad \quad  X^n(t=0,x)=x. 
\end{equation} 
We also denote by $X$ the a.e flow of $u$. 

Let $Q^n$ be the solution of \eqref{NF} with $F^n$ replaced by $|F^n|$ and 
taking the same initial data as $N^n_2$ at $t=0$. Hence, 
\begin{equation} \label{NF-renor-flow} 
\frac{d  [\beta_M(Q^n) (t,X^n(t,x) ) ] }{dt}  
 = \beta_M'(Q^n)  | F^n (t,X^n(t,x) ) |    
\end{equation}  
 where the equation holds in the sense of distribution. 

Passing to the limit weakly in \eqref{NF-renor-flow}, we get 

\begin{equation} \label{NF-renor-flow-lim}
\frac{d \overline{[\beta_M(Q^n) (t,X^n(t,x) ) ]} }{dt}   =  
 \overline {\beta_M'(Q^n)   |F^n (t,X^n(t,x) )|  }   
\end{equation}

From the stability of the notion of   a.e flow with respect to 
the  weak limit 
of $u^n$ to $u$, we know  that $ X^n(t,x) $ converges to  $ X(t,x) $ in 
$L^1_{loc}$ and also that 
$ (X^n (t)^{-1}) (x) $ converges to  $  (X  (t)^{-1}) (x)   $ in 
$L^1_{loc}$.
 This allows us to get the following equality of the  weak limits
\begin{equation} 
 \overline{[\beta_M(Q^n) (t,X^n(t,x) ) ]}   =  
 \overline {  [\beta_M(Q^n) (t,X(t,x) ) ] }.    
\end{equation}  
 
Now, sending $M$ to infinity in \eqref{NF-renor-flow-lim}, we deduce that 

\begin{equation} \label{NF-lim}
\frac{d [Q  (t,X (t,x) ) ] }{dt}    =   F
\end{equation}  

where $Q =  \lim_{M \to \infty}  \overline {  [\beta_M(Q^n) ] }   $   is the 
Chacon limit of $Q^n$ and 
$F = \lim_{M \to \infty} \overline {\beta_M'(Q^n)   |F^n (t,X^n(t,x) )|  }    $.
 It is easy to see that $Q \in L^\infty (0,T ; L^1 (\Omega))$ and 
that $F  \in \M ((0,T) \times \Omega) $. Integrating in $t$, we deduce that 
a.e in $x\in \Omega$, we have

\begin{equation} \label{Q=}
 Q  (t,X (t,x) )     = Q(0,x) +    \int_0^t  F (s,X (s,x) ) \, ds 
\end{equation}  
for a.e $t\in (0,T)$.  Due to the fact that $F $ is nonnegative, we deduce that 
$ Q  (t,X (t,x) )    $ is increasing in time. 
We define  the normalizing factor $N$  by the following 

\begin{equation} \label{Q=1}
N  (t,X (t,x) )  =Q  (T_0,X (T_0,x) )   = Q(0,x) +    \int_0^{T_0}  F (s,X (s,x) ) \, ds 
\end{equation}  
for  $t\in (0,T_0)$
where $T_0< T$ is a fixed time. In the sequel, we will denote $T=T_0$
and will not make the distinction between these two times. 
Notice that  
$N$ is  constant  along the characteristics of $u$,  that $N$ is in $ L^\infty (0,T ; L^1 (\Omega)) $ 
and that $N  (t,X (t,x) )   $ is in  $   L^1 ( \Omega; L^\infty (0,T ) )    $. 
Moreover it is bounded from below by $1$. Hence, it solves 
$$(\partial_t + u.\nabla) \frac1N = 0 .  $$
Also,  the following  two inequalities  hold

 \begin{equation} 
\overline{ \beta_M(N_2^n) }   \leq  \overline{ \beta_M(Q^n)}  \leq Q \leq N   
\end{equation}  
 and hence the weak limit of $N_2^n$  which is equal to 
  the Chacon limit of $N_2^n$ is  bounded by $N$.  The fact that 
 the weak limit of $N_2^n$   is equal to 
 its   Chacon limit comes from the fact that  the sequence  $N_2^n$ is equiintegrable. This is 
a  simple consequence of the dissipation of the free energy and the 
weighted Sobolev inequality \eqref{WSI-eq}.
Indeed,  from \eqref{WSI-eq}, we deduce that $\sqrt{ \psi^n \over \psi_\infty }$ is 
bounded in  $ L^2( (0,T) \times \Omega ;  L^p( \psi_\infty dR)   ) $ on the other 
hand from the conservation of mass, we know that 
 $\sqrt{ \psi^n \over \psi_\infty }$ is 
bounded in  $ L^\infty ( (0,T) \times \Omega ;  L^2( \psi_\infty dR)   ) $. 
Interpolating between these two bounds, we easily deduce that  
$\sqrt{ \psi^n \over \psi_\infty }$ is 
bounded in  $ L^r  ( (0,T) \times \Omega \times B, dt dx    \psi_\infty dR    ) $ for 
some $r> 2$ and hence  
$N_2^n$ is equiintegrable. We also get that 
$N^n_1 $ is equiintegrable and hence  $\theta$ which is 
 the Chacon limit of $N^n_1 $  is equal to  the weak limit of $N^n_1 $.


\subsection{The term  $\overline {  \nabla u^n : {\tau (\tilde \psi^n) } }^{\de,\ka} $  }

In this subsection, we will prove that $ \overline {  \nabla u^n : {\tau (\tilde \psi^n) } }^{\de,\ka}  
   = \nabla u : \tau + \int_B \beta_{ij} \frac{R_i R_j}{1-|R^2|} .  $ This will follow 
from the following two lemmas

\begin{lem}  \label{L1}
\begin{equation}  \label{eq-L1}
\overline {  \frac{ \nabla u^n : {\tau (\tilde \psi^n) } }{ (1 + \delta N_2^n ) (1 + \ka N_1^n)^2    } } 
 =   \int_B      z^{\de,\ka}        \frac{R_i R_j}{1-|R^2|} 
\end{equation}  
where $z^{\de,\ka} = \overline {  \frac{  {  \tilde \psi^n } \,   \nabla u^n   }{ (1 + \delta N_2^n ) (1 + \ka N_1^n)^2    } }       $ 
\end{lem}

\begin{lem}\label{L2}
  $z^{\de,\ka}$  converges strongly to $ \overline { \nabla u^n   \tilde \psi^n } = 
   \nabla u    \tilde \psi  + \beta    $  in $ L^1( (0,T) \times \Omega \times B; 
dt \, dx \, \frac{dR}{1-|R|}  )  $ when $\delta$ goes to zero and then $\ka$ goes to zero. 
\end{lem}
 Denoting $\tau^{n,\de,\ka} = \frac{ {\tau (\tilde \psi^n) } }{ (1 + \delta N_2^n ) (1 + \ka N_1^n)^2    }   $, we get that 
\begin{cor} \label{u-t-n} 
\begin{equation}  
 \overline { \nabla u^n  : \tau (\tilde \psi^n)}^{\de, \ka} 
 = \lim_{\ka \to 0} \lim_{\de \to 0} 
 \overline { \nabla u^n  : \tau^{ n,\de,\ka}}   = \nabla u : \tau(\psi) 
 + \int_B \be_{ij}  \frac{R_i R_j}{1-|R^2|} dR. 
\end{equation}  
\end{cor}

{\bf Proof of Lemma \ref{L1}.} 
The proof of \eqref{eq-L1} follows from the fact that 
$ z^{n, \de,\ka}   =  \frac{ \nabla u^n    \tilde \psi^n  } { \psi_\infty    (1 + \delta N_2^n ) (1 + \ka N_1^n)^2    }  $
is equi-integrable in $ L^1( (0,T) \times \Omega \times B; 
dt \, dx \, \frac{ \psi_\infty dR}{1-|R|}  )  $  for $\de, \ka$ fixed. 
Indeed, consider the real valued  function   $\Phi(x) = x \log (1+x)  + 1 $. It is enough 
to prove that 
$\Phi (z^{n, \de,\ka}  )  =   \Phi \left( \frac{ \nabla u^n    \tilde \psi^n  } { \psi_\infty    (1 + \delta N_2^n ) (1 + \ka N_1^n)^2    }  \right)  $
is bounded  in $ X =  L^1( (0,T) \times \Omega \times B; 
dt \, dx \, \frac{ \psi_\infty dR}{1-|R|}  )  $. To simplify notation, we denote 
$N^n = (1 + \delta N_2^n ) (1 + \ka N_1^n)^2   $. Hence, it is enough to bound

\begin{equation}  
 \label{2-terms}
   \frac{ \nabla u^n}{N^n}  \left[ 
    \frac{ \tilde \psi^n  } { \psi_\infty }   
 \log \left(   \frac{ \tilde \psi^n  } { \psi_\infty }  \right)    
+  \frac{ \tilde \psi^n  }{ \psi_\infty }  \log \left(    \frac{ \nabla u^n}{N^n}   \right) 
\right]
\end{equation}       
in $X$ (see definition above).

To bound the first term appearing in \eqref{2-terms} we use  the  Hardy type inequality 
\eqref{Hardy-inter} to get that 

\begin{align} \label{1st-terms}
&  \frac{ \nabla u^n}{N^n}  \int_B  
     \tilde \psi^n  
 \log \left(   \frac{ \tilde \psi^n  } { \psi_\infty }  \right)    
\frac{1}{1-|R|} \ dR  \\ 
& \quad \les \frac{ \nabla u^n}{N^n}  \left[  \int_B  
     \tilde \psi^n     
 \log \left(   \frac{ \tilde \psi^n  } { \psi_\infty }  \right)        \right]^{1/2}     
\left[   \int_B   \psi_\infty \left|  \nabla \Big( \sqrt{ \frac{ \tilde \psi^n  } { \psi_\infty }  } 
  \log^{1/2}  \frac{ \tilde \psi^n  } { \psi_\infty }      \Big)  \right|^2       \right]^{1/2} \\ 
& \quad \les | \nabla u^n |^2  + \frac{ 1}{ (N^n)^2} \left[  \int_B  
     \tilde \psi^n     
 \log \left(   \frac{ \tilde \psi^n  } { \psi_\infty }  \right)        \right]      
\left[   \int_B   \psi_\infty \left|  \nabla \Big( \sqrt{ \frac{ \tilde \psi^n  } { \psi_\infty }  } 
  \log^{1/2}  \frac{ \tilde \psi^n  } { \psi_\infty }      \Big)  \right|^2       \right] 
\end{align}    
 and using the a priori bound \eqref{log2-b}, we see that the last term 
in uniformly bounded in $L^1( (0,T) \times \Omega  ) $.

To bound the second term in \eqref{2-terms}, we first  use the inequality 
$ x\, y  \leq  C ( x^2 \log^2(x)  + \frac{y^2}{\log^2 y }  )    $ for $x,y \geq 2 $ and 
then apply Jensen inequality.  
Hence, 
 
\begin{align} \label{2nd-terms}
&  \frac{ \nabla u^n}{N^n}  \log \left(    \frac{ \nabla u^n}{N^n}   \right)  \int_B  
    \frac{ \tilde \psi^n  } { \psi_\infty }   
\frac{\psi_\infty}{1-|R|} \ dR  \\ 
& \quad \les  | \nabla u^n |^2  + \frac{ 1}{ (N^n)^2} \left[  \int_B  
     \frac{\tilde \psi^n  }{ \psi_\infty }   
\frac{\psi_\infty}{1-|R|} \ dR  \right]^2   
\log^2 \left[  \int_B  
     \frac{\tilde \psi^n  }{ \psi_\infty }   
\frac{\psi_\infty}{1-|R|} \ dR  \right]  \\
& \quad \les  | \nabla u^n |^2  + \frac{ 1}{ (N^n)^2} \left[  \int_B  
     \frac{\tilde \psi^n  }{ \psi_\infty }  \log \left(   \frac{ \tilde \psi^n  } { \psi_\infty }  \right)      
\frac{\psi_\infty}{1-|R|} \ dR  \right]^2  
\end{align}    
and the last term can be bounded as in \eqref{1st-terms}. We notice here 
that the last inequality implies in particular that 
$|\tau^{n,\de,\ka}|^2$ is equi-integrable in $L^1$ for fixed $\de$ and $\ka$.  
This is actually a very important fact that will be used again later.

{\bf Proof of Lemma \ref{L2}.} 
To prove this lemma, we use dominated convergence and monotone convergence. Indeed, 
$ |  z^{n, \de,\ka}    |  $  is decreasing in $\delta, \ka$, namely
for $ 0< \de \leq  \de'$ 
and $0< \ka \leq \ka'$, we have 

\begin{equation}
 |  z^{n, \de',\ka'}    | \leq  |  z^{n, \de,\ka}    |  \leq  | \nabla u^n  { (\tilde \psi^n) }     |.
\end{equation}
Passing to the limit weakly in $n$, we deduce that 
\begin{equation}
\overline{ |  z^{n, \de',\ka'}    | }  \leq  
 \overline{ |  z^{n, \de,\ka}    | }   \leq  
\overline{ | \nabla u^n  { (\tilde \psi^n) }     | } 
\end{equation}
 and by monotone convergence, we deduce that 
$ G  =  \overline{ |  z^{n, \de,\ka}    | }^{\de,\ka} \in X    $ and
 that for all   $ 0< \de $ 
and $0< \ka $, we have 
$|  z^{\de,\ka}    | \leq G  $. Moreover, we have 

\begin{equation}
|  z^{\de,\ka}  -   z^{\de',\ka'}   | \leq \left|  \overline{ |  z^{n, \de,\ka}    | }  
-   \overline{ |  z^{n, \de',\ka'}    | }      \right|.  
\end{equation}   
 Hence, there exists $g\in X$ such that $ z^{\de,\ka} $ converges strongly 
to $g$ in $X$. 
Now, we would like to prove that the limit $g$ is equal to $ \overline { \nabla u^n   \tilde \psi^n } $.
This follows from the fact that $  \nabla u^n   \tilde \psi^n  $ is equi-integrable in 
$L^1((0,T) \times \Omega \times B; 
dt \, dx \, dR ) $  (without the weight). Indeed, denoting 
$\Phi(x) =|x| \log^{1/2} (1+|x|)   $, we have 
\begin{align}
\Phi( | \nabla u^n   \tilde \psi^n   |) & \les |\nabla u^n   \tilde \psi^n   | ( \log (1+ \tilde \psi^n) 
  + \log(1 +  | \nabla u^n |   )  ) \\
& \les  \tilde \psi^n [ |\nabla u^n   |^2 + \log(\tilde \psi^n)    ]
\end{align}   
which is clearly bounded in $L^1(dt\, dx\, dR)$.

\subsection{Identification of  $\int_B \beta_{ij}  \frac{R_i R_j}{1-|R^2|} dR $.}
In this subsection, we give a relation between $\be$ and some  defect 
measure related to the lack of strong convergence of $\nabla u^n$ in $L^2$. 
To state the  main proposition  of this subsection, we introduce few notations. 
Let $u^n = v^n + w^n$ where $v^n$ and $w^n $solve 
\begin{align} \label{v-n}
 \left\{     \begin{array}{l} 
  \partial_t v^n - \Delta v^n + \nabla p_1^n = \nabla. \tau^n  \\ 
   v^n (t=0) = 0  
  \end{array}\right. 
\end{align}   
\begin{align}\label{w-n}
 \left\{     \begin{array}{l} 
  \partial_t w^n - \Delta w^n + \nabla p_2^n =  - u^n.\nabla u^n   \\ 
   v^n (t=0) = u^n(t=0). 
  \end{array}\right. 
\end{align}   
We further split $w^n$ into  $w^n_1 + w^n_2$ where $w^n_1$ is 
the solution with zero initial data and $w^n_2$ is 
the solution with zero right hand side. 

In the rest of this subsection we will use   $\De$  to   denote 
 $\de,\ka $. 
We  define $   v^{n,\De} =   v^{n,\de,\ka}$ the solution of 
\begin{align} \label{v-n-d-k}
 \left\{     \begin{array}{l} 
  \partial_t v^{n,\De} - \Delta v^{n,\De} + \nabla p_1^{n,\De} = \nabla. \tau^{n,\De}  \\ 
   v^{n,\De} (t=0) = 0   
  \end{array}\right.  
\end{align}    
Extracting a subsequence, we assume that 
$( \tau^{n,\De} , \nabla   v^{n,\De}, \nabla v^n, \nabla w^n)  $
 converges weakly  in $L^2$ to some $( \tau^{\De} ,   \nabla   v^{\De}, \nabla v,\nabla w)
  $ and that 
\begin{equation} 
\overline{| \nabla   v^{n,\De}  |^2 } = |\nabla   v^{\De}    |^2  + \mu^{\De}
\end{equation}   
for some defect measure $\mu^{\De} \in \M ( (0,T)\times \Omega ) $. 
We also denote $\mu$ the limit of $\mu^{\De}$  when $\de$ and then
$\ka$ go to zero (extracting a subsequence), namely 
\begin{equation} 
 \mu =  \lim_{\ka \to 0}  \lim_{\de \to 0}    \mu^{\De} =  \lim_{\De \to 0}    \mu^{\De} . 
\end{equation}

\begin{prop} \label{mu-beta}
We have 
\begin{equation} 
  \mu  = -  \int_B \beta_{ij}  \frac{R_i R_j}{1-|R^2|} dR. 
\end{equation}  
\end{prop}

{\bf Proof of Proposition \ref{mu-beta}.}

We   introduce the following weak limits 
\begin{align} 
\overline{ \tau^{n,\De} : \nabla v^{n,\De}  } &= W^{\De\De} \\ 
\overline{ \tau^{n,\De} : \nabla v^{n}  } &= W^\De. 
\end{align}  
 
{\bf Step 1:}  
First, we would like to prove that $W^{\De\De} $ and $ W^\De $ have the same 
limit $W$  when $\De$ goes to zero and that this limit is in $L^1$. 
To prove  this, we  introduce for $M > 0$, the following weak limits   
\begin{align} 
\overline{ \tau^{n,\De} 1_{ |\tau^{n,\De} | \leq M  } : \nabla v^{n,\De}  } &= W^{\De\De}_M  \\ 
\overline{ \tau^{n,\De} 1_{ |\tau^{n,\De} | >  M  } : \nabla v^{n,\De}  } &=  W^{\De\De} -    W^{\De\De}_M  \\
\overline{ \tau^{n,\De} 1_{ |\tau^{n,\De} | \leq M  }  : \nabla v^{n}  } &= W^\De_M
\end{align}  
and 
\begin{align} 
\overline{  |\tau^{n,\De} 1_{ |\tau^{n,\De} | \leq M }  |^2   }  = G^{\De}_M   \quad \quad 
  \overline{  |\tau^{n,\De}  |^2 } = G^{\De} .     
\end{align}  
Since for a fixed $\De$, $ |\tau^{n,\De}  |^2  $  is equi-integrable, we deduce that 
$ G^{\De}_M  $ converges to $G^{\De}$ in $L^1$ when $M$ goes to infinity and is monotone 
in $M$. Also, by monotone convergence, we deduce that there exists $G \in L^1$ such 
that $G^{\De}$    converges to $G$ in $L^1$ when $\De$ goes to zero. Actually, 
$G$ is the weak limit of $|\tau^n|^2 $ in the sense of Chacon. 

Let us fix $\eps > 0$. We choose 
$\De_0$ and $M_0$ such that for $\De < \De_0 $ and $M > M_0$, we have 
$\| G-G^\De  \|_{L^1} + \| G-G^\De_M  \|_{L^1} \leq \eps    $. We have 

\begin{align} 
\overline{ |\tau^{n,\De} |^2 } & = \overline{  |\tau^{n,\De} 1_{ |\tau^{n,\De} | \leq M }  |^2   }  
 + \overline{  |\tau^{n,\De} 1_{ |\tau^{n,\De} | >  M }  |^2   } \\ 
 & =   G^{\De}_M  + \quad (G^\De- G^{\De}_M  ). 
\end{align}  
Hence, we deduce that for $\De < \De_0 $ and $M > M_0$, we have  for all $n$, 
$\| \tau^{n,\De} 1_{ |\tau^{n,\De} | >  M }     \|_{L^2}^2  \leq \eps  $ and hence, by 
Cauchy-Schwarz  we
deduce that $  \|  W^{\De\De} -    W^{\De\De}_M \|_{L^1} \leq C \sqrt{\eps}   $ and 
that $  \|  W^{\De} -    W^{\De}_M \|_{L^1} \leq C \sqrt{\eps}   $.  
Hence to prove that $\lim_{\De} W^{\De\De} = \lim_{\De} W^{\De}  $, it is enough 
to prove it for the $M$ approximation, namely that 
\begin{equation}  \label{M-lim} 
 \lim_{\De} W^{\De\De}_M = \lim_{\De} W^{\De}_M  . 
\end{equation}  

To prove \eqref{M-lim}, we first notice that 
$\tau^{n,\De}  - \tau^n  $ goes to zero in $L^p$ for $p< 2$ 
 when $\De$ goes to zero uniformly in 
$n$. Then, by parabolic regularity of the Stokes system, we deduce that  
$ \| \nabla v^{n,\De} -  \nabla v^{n}       \|_{L^p( (0,T) \times \Omega    )  }  $ goes to zero when 
$\De$ goes to zero uniformly in $n$ for $p< 2$. Hence, \eqref{M-lim} holds. 

{\bf Step 2:} In this second step, we will compare the local  energy identity of  the weak limit of 
\eqref{v-n-d-k}  with the weak limit of the local energy identity of \eqref{v-n-d-k}. 

On one hand, passing to the limit in \eqref{v-n-d-k} and multiplying by $v^\De$, we deduce that 
\begin{equation} \label{v-n-d-k-limit}
   \partial_t \frac{|v^{\De}|^2}2  - \Delta \frac{|v^{\De}|^2}2  +  
|\nabla v^{\De} |^2     + \div ( p_1^{\De} v^{\De} )   = 
 \div ( v^\De .  \tau^{\De}) -  \nabla v^\De : \tau^\De   
\end{equation}  

On the other hand, reversing the order, we get 

\begin{equation} \label{v-n-d-k-limit2}
   \partial_t \frac{|v^{\De}|^2}2  - \Delta \frac{|v^{\De}|^2}2  +  
|\nabla v^{\De} |^2  + \mu_\De    + \div ( p_1^{\De} v^{\De} )   = 
 \div ( v^\De .  \tau^{\De}) -   W^{\De \De}. 
\end{equation}  
For a justification of these two calculations, we refer to \cite{LM00cam}. 
Comparing \eqref{v-n-d-k-limit} and \eqref{v-n-d-k-limit2}, we deduce that 
$ W^{\De\De} = \nabla v^\De : \tau^\De   - \mu_\De  $. 
We would like now to send  $\De$ to zero. 
 
First, it is clear that $\tau^\De$ converges strongly to $\tau$ in $L^2$ when
$\De$ goes to zero. Hence, $\nabla v^\De$ also converges to  $\nabla
v$ in $L^2$. Besides, from the energy estimate, we recall that 
$u^n $ is bounded in $L^\infty((0,T); L^2(\Omega)) \cap  L^2((0,T); \dot H^1(\Omega)) $ 
and hence by Sobolev embeddings that $u^n  $ is bounded in $L^{\frac{2(D+2)}{D}} ((0,T)\times \Omega )  $
and that $u^n \nabla u^n $ is bounded in $L^{\frac{D+2}{D+1}} ((0,T)\times \Omega )  $. 
By parabolic regularity of the Stokes operator applied to 
\eqref{w-n} with zero initial data, we deduce that $\nabla w^n_1$ 
is bounded in  $ L^{\frac{D+2}{D+1}} ((0,T); W^{1,\frac{D+2}{D+1} }\Omega )  $ and 
that $ \partial_t w^n_1$ is bounded in $L^{\frac{D+2}{D+1}} ((0,T)\times \Omega )  $. 
 Since $\tau^n $ is bounded in $L^2$, we deduce  from \eqref{v-n} that 
$\nabla v^n$ is also bounded in $L^2((0,T)\times \Omega ) $ and hence 
$\nabla w^n$ is also bounded in $L^2((0,T)\times \Omega ) $. 
Moreover, it is clear that $\nabla w^n_2$ is compact in  $L^2((0,T)\times \Omega ) $
and hence $\nabla w^n_1$ is also bounded in $L^2$ and from the previous 
bounds on $ \nabla w^n_1$, we deduce that 
$ \nabla w^n_1$ is compact in  $ L^p((0,T)\times \Omega ) $ for $p<2$. 
Hence, we deduce that 
$  \overline { \nabla w^n  : \tau (\tilde \psi^n)}^{\de, \ka}  = \nabla w : \tau(\psi)   $ 
(where we have used that $\tau^{n,\De}$ is equi-integrable for each fixed $\De$)
and from Corollary \ref{u-t-n}  that 
$  \lim_{\De }  W^{\De \De}  =  \lim_{\De }  W^\De = 
 \overline { \nabla v^n  : \tau (\tilde \psi^n)}^{\de, \ka}  = \nabla v : \tau(\psi)   
 + \int_B \be_{ij}  \frac{R_i R_j}{1-|R^2|} dR. $ 
Finally, we deduce that 
$ \mu = \lim_{\De \to 0 } \mu^\De  = -  \int_B \be_{ij}  \frac{R_i R_j}{1-|R^2|} dR. $


\subsection{Gronwall along the characteristics} \label{Gron-sec}

Taking  the difference between \eqref{psi-n-delta} and
\eqref{g2log-eq} and dividing by $N^4$, we  get  (to be more precise, we have to take 
the difference between \eqref{psi-n-delta-De}  and \eqref{g2log-eq}, divide by $N^4$ and 
then send $\De$ to zero):

\begin{align}
\nonumber & 
 (\partial_t + u . \nabla) \frac{ \eta  }{N^4}   \\ 
     & \quad =   
  \frac1{N^4}  \left[ \overline {  \nabla u^n : {\tau (\tilde \psi^n) }
    }^{\de,\ka}  +  \nabla u : \tau\left( \psi_\infty \Big(    g^2 - 2g \overline{\pp
     \Th'(\pp) }  \Big)  \right)   
                        \right]\nonumber   \\ 
 \label{Gron-eq}   & \quad  \quad  -   \frac1{N^4}  \int_B \psi_\infty  \left[ 4 \overline { 
       \left|    \nabla_R  f^n      \right|^2  }^{\de,\ka}   
 -2 |\nabla_R g|^2  + \overline{ \frac{|\nabla_R f^n|^2 (\log^{1/2} + \log^{-3/2})(\pp) }{f^n}} 2g                       \right] \\ 
  &\quad  \quad  -   \frac{2ak}{N^4}  \int_B  \left(  \overline {\nabla  u^n   
\Big(\log( \frac{\tilde \psi^n}{\psi_\infty} ) + 1  \Big)
  }^{\de,\ka}   -  
  \Big(   \o{ 2 \Th'(\pp) \nabla_i u^n_j   }  \Big) g   \right)
  \frac{R_iR_j}{1-|R|^2} \psi_\infty                 \nonumber   \\  
 &\quad  \quad  -   \frac2{N^4} \int \psi_\infty  \left[ \gamma_{ij} R_j \nabla g
     - \nabla_R \U  R : (\gamma -\gamma') g           \right]\nonumber
  \\
  & \quad  =  -  \sum_{i=1}^4 A_i \nonumber  
\end{align}    
where we  denote the 4 terms appearing on the right hand side by $A_i,
1\leq i \leq 4$ and we also 
 denote $\eta = \overline { N^n_1}^{ \de,\ka}  - \int_B g^2 \psi_\infty = \int_B 
  [  \o{ (g^n)^2 }  -g^2 ] \psi_\infty dR    $. 
 It measures  the lack of strong convergence of $g^n  $
to $g$ in $L^2(dtdx \psi_\infty dR)$. 
 Notice that by the choice of the normalizing factor $N$, the defect measure 
$\frac{\eta}N$ is in $L^\infty$. 
  
First, we prove  that 
 $A_2$ is  nonnegative, namely 
we have the following lemma

\begin{lem}  \label{def-A2}
We have 
\begin{equation}
A_2 \geq  \frac{c}{N^4} \int_B \psi_\infty  \overline { 
       \left|    \nabla_R  (f^n -f)       \right|^2  }^{\de,\ka}
 =  \frac{c}{N^4} \int_B \psi_\infty  \varpi 
\end{equation}
for some  constant $c$. 
\end{lem}

For the proof, we 
rewrite $ |\nabla_R g|^2  $ as

 \begin{align} 
 |\nabla_R g|^2  & =   \left| \o{\nabla_R f^n (\log^{1/2} (f^n)^2 +
 \log^{-1/2} (f^n)^2   )  } \right|^2  \\ 
\label{ng2} &   =  \left| \o{\nabla_R f^n (\log^{1/2} (f^n)^2)}  \right|^2
  +   \left| \o{\nabla_R f^n (\log^{-1/2} (f^n)^2)}  \right|^2 \\
& \quad  + 2  \o{\nabla_R f^n (\log^{1/2} (f^n)^2)}  \cdot  
 \o{\nabla_R f^n (\log^{-1/2} (f^n)^2)}.     
\end{align}

\def\a{\boldsymbol\alpha }
\def\b{\boldsymbol\beta }
\def\c{\boldsymbol\gamma }

Hence, we deduce that 
\begin{equation}
A_2 =    \frac1{N^4}  \int_B \psi_\infty   (\a+ \b + \c)  
\end{equation} 
where $\a,\b$ and $\c$ are given by  
 
\begin{equation}
\frac{\a}2 = \overline{ \frac{|\nabla_R f^n|^2 (\log^{1/2} (\pp)  ) }{f^n}}    
 \o{ f^n  \log^{1/2} (\pp)   }
  - \o{  (\nabla f^n) \log^{1/2} (\pp)  }^2   
\end{equation}

\begin{equation}
\frac{\b}2 = \overline{ \frac{|\nabla_R f^n|^2 (\log^{-3/2} (\pp)  ) }{f^n}}    \o{ f^n  \log^{1/2} (\pp)   }
  - \o{  (\nabla f^n) \log^{-1/2} (\pp)  }^2   
\end{equation}

\begin{equation}
\frac{\c}2 =  2 \overline{ |  \nabla f^n |^2     }^{\de,\ka}    
  -  2   \o{  (\nabla f^n) \log^{1/2} (\pp)  }  \quad  \o{  (\nabla f^n) \log^{-1/2} (\pp)  }
\end{equation}

We introduce the Young measure  $\nu_{t,x,R}(\Lambda, \lambda)$ associated to the sequence 
$(\nabla f^n, f^n)$ where $\Lambda \in \R^D$  and $\lambda \in \R$. 
Hence, the defect measure  $ \overline { 
       \left|    \nabla_R  (f^n -f)       \right|^2  }^{\de,\ka}  $
satisfies : 

\begin{align}
  \overline { 
       \left|    \nabla_R  (f^n -f)       \right|^2  }  \geq  \overline { 
       \left|    \nabla_R  (f^n -f)       \right|^2  }^{\de,\ka}
  &\geq 
        \int   |\Lambda - \int \Lambda'  \nu_{t,x,R} (\Lambda',\lambda')   |^2   \nu_{t,x,R} (\Lambda,\lambda) \\ 
     &   =  \frac12     \int  \int   |\Lambda - \Lambda' |^2    \nu_{t,x,R} (\Lambda',\lambda')     \nu_{t,x,R} (\Lambda,\lambda)  
\end{align}
Indeed, it is easy to see that  $\overline { 
       \left|    \nabla_R  (f^n -f)       \right|^2  }^{\de,\ka} $ is
bounded from above by the weak limit and from below by the Chacon
limit of  $ \left|    \nabla_R  (f^n -f)       \right|^2  $.  
  In the sequel, we will drop the $t,x$ and $R$ dependence of $\nu$ and will denote 
  $\nu' = \nu(\Lambda',\lambda')   $  and  $\nu = \nu(\Lambda,\lambda)   $.   
  Besides, $\a,\b$ and $\c$ satisfy  
  \begin{align} 
\a \geq    \int  \int    A(   \Lambda ,\lambda, \Lambda',\lambda'   )    \nu (\Lambda',\lambda')     \nu (\Lambda,\lambda)  
\end{align} 
and the same for $\b$ and $\c$ with $A$ replaced by $B$ or $C$ where $A, B$ and $C$ are 
given by 

  \begin{align}
A & = \frac{|\Lambda|^2 \log^{1/2} (\lambda^2)     }{\lambda}   \lambda' \log^{1/2} (\lambda')^2
  +     \frac{|\Lambda'|^2 \log^{1/2} (\lambda')^2     }{\lambda'}   \lambda \log^{1/2} (\lambda^2)
   - 2 \Lambda .\Lambda'  \log^{1/2} (\lambda^2)   \log^{1/2} (\lambda')^2  \\
   B& =  \frac{|\Lambda|^2 \log^{-3/2} (\lambda^2)     }{\lambda}   \lambda' \log^{1/2} (\lambda')^2
  +    \frac{|\Lambda'|^2 \log^{-3/2} (\lambda')^2     }{\lambda'}   \lambda \log^{1/2} (\lambda^2)
   - 2\Lambda .\Lambda'  \log^{-1/2} (\lambda^2)   \log^{-1/2} (\lambda')^2  \\
   C & =   2 |\Lambda|^2   +  2  | \Lambda' |^2
   - 2 \Lambda .\Lambda'    \Big(  \log^{1/2} (\lambda^2)   \log^{-1/2} (\lambda'^2)   +     \log^{-1/2} (\lambda^2)   \log^{1/2} (\lambda')^2       \Big) . 
\end{align} 
  To prove lemma \ref{def-A2}, it is enough to show that 
  $A+B+C \geq  \frac{c}2    |\Lambda - \Lambda'|^2$. 
  First, we rewrite $A+B+C $ as 
  \begin{align}
   A+B+C &= |\Lambda|^2 B_1  +   |\Lambda'|^2 B_2  -2 \Lambda. \Lambda' B_3  \\ 
               & =  |  \Lambda - \Lambda'  |^2  +  |\Lambda|^2 (B_1-1)   +  
                  |\Lambda'|^2   (B_2-1)   -2 \Lambda. \Lambda'   (B_3 -1)    
  \end{align} 
  where $B_1, B_2$  and $B_3$
 are given by 
      \begin{align}
B_1 & = \frac{\log^{1/2} (\lambda^2)     }{\lambda}   \lambda' \log^{1/2} (\lambda')^2
  + \frac{  \lambda' \log^{1/2} (\lambda')^2   }{\lambda \log^{3/2} (\lambda^2)    }  + 2 \\ 
B_2 & = \frac{\log^{1/2} (\lambda'^2)     }{\lambda'}   \lambda \log^{1/2} (\lambda^2)
  + \frac{  \lambda \log^{1/2} (\lambda^2)    }{\lambda' \log^{3/2} (\lambda'^2)    }  + 2 \\   
  B_3 & =   \log^{1/2} (\lambda^2)  \log^{1/2} (\lambda'^2)  + \frac1{  \log^{1/2} (\lambda^2)  \log^{1/2} (\lambda'^2)  }    + \frac{  \log^{1/2} (\lambda^2)  } {  \log^{1/2} (\lambda'^2)   }  
 + \frac{  \log^{1/2} (\lambda'^2)  } {  \log^{1/2} (\lambda^2)   }   
\end{align}

  Actually, we will prove that if $a$ is chosen big enough then $ (B_1 - 1)(B_2-1) \geq 
   (B_3-1)^2 $ from which we deduce that $A+B+C \geq | \Lambda -
  \Lambda'  |^2$ and the lemma would follow. 
   Indeed, after simple calculations, we get  
    
      \begin{align*}
& (B_1 - 1)(B_2-1)  -  
   (B_3-1)^2     = \\
    & \quad \quad  \log^{1/2} (\lambda^2)  \log^{1/2} (\lambda'^2) \left[  \frac{\lambda}{\lambda'} 
    + \frac{\lambda'}{\lambda}     + 2 -   2  \frac{  \log^{1/2} (\lambda^2)  } {  \log^{1/2} (\lambda'^2)   }  
  -2   \frac{  \log^{1/2} (\lambda'^2)  } {  \log^{1/2} (\lambda^2)   }        \right]  \\
  &  \quad \quad +   2 \left[   \frac{  \log^{1/2} (\lambda^2)  } {  \log^{1/2} (\lambda'^2)   }  
 + \frac{  \log^{1/2} (\lambda'^2)  } {  \log^{1/2} (\lambda^2)   }     - 2    \right]  \\ 
 &\quad \quad +  \frac1{  \log^{1/2} (\lambda^2)  \log^{1/2} (\lambda'^2)  }    \left[ \frac{   \lambda  \log(\lambda^2)  } {    \lambda'  \log (\lambda'^2)   }  
 + \frac{   \lambda'  \log (\lambda'^2)  } {   \lambda  \log (\lambda^2)   }   + 
 2   - 2    \frac{  \log^{1/2} (\lambda^2)  } {  \log^{1/2} (\lambda'^2)   }  
 -2  \frac{  \log^{1/2} (\lambda'^2)  } {  \log^{1/2} (\lambda^2)   }     \right]
\end{align*}  
We will prove that the three terms appearing inside the brackets are nonnegative. 
This is obvious for the second one since it is of the form $ x +\frac1x -2 $ for some $x>0$. 
We recall that since $(f^n)^2 \geq a$, we get that $\lambda \geq \sqrt{a}$ on the 
support of $\nu$.  For the first bracket, we assume that $\lambda' \geq \lambda$ and 
write $\lambda' = \lambda (1 +\eps) $. Hence, the term in the first
bracket   is given by 

\begin{equation}  \label{1-eps}
   1+\eps + \frac1{1+\eps} + 2 -  2  \sqrt{ 1 +  \frac{\log(1+\eps) }{  \log \lambda  }  } 
    -2 \frac1{   \sqrt{ 1 +  \frac{\log(1+\eps) }{  \log \lambda  }  }  } 
    \end{equation}
  and one can check easily that if $\lambda \geq \sqrt{a}$ is big enough then 
  \eqref{1-eps}   is nonnegative.  The same argument can be used for the third 
bracket.  This end the proof of lemma \ref{def-A2}. 
  \vspace{.2cm}
  
  To bound  $A_1$, we first observe that 
\begin{equation}  
g^2 - 2g \overline{\pp
     \Th'(\pp) }  = -  \o{f^n \log^{1/2} (f^n)^2}   \  \o{f^n \log^{-1/2} (f^n)^2}
    \end{equation}
  and hence, 

 \begin{align} 
A_1  & =  - \frac1{N^4} \left[  \overline {  \nabla u^n : {\tau (\tilde \psi^n) }
    }^{\de,\ka}  +  \nabla u : \tau\left( \psi_\infty \Big(    g^2 - 2g \overline{\pp
     \Th'(\pp) }  \Big)  \right)  \right]  \\ 
   &   = - \frac1{N^4} \left[   \overline {  \nabla u^n : {\tau (\tilde \psi^n - \psi ) }
    }^{\de,\ka}  + \nabla u : \tau ( \psi - \psi_\infty ( \o{f^n \log^{1/2} (f^n)^2}   \  \o{f^n \log^{-1/2} (f^n)^2}   ) ) \right] \\ 
 &  =   \frac1{N^4} \left[   \mu - \nabla u : \tau ( \psi - \psi_\infty ( \o{f^n \log^{1/2} (f^n)^2}   \  \o{f^n \log^{-1/2} (f^n)^2}   ) ) \right]   
\end{align}     
By convexity, it is clear that 
$\o{ (f^n -f )^2   }  = \o{(f^n)^2 } - f^2 \geq  \o{(f^n)^2 }  
 - \o{f^n \log^{1/2} (f^n)^2}   \  \o{f^n \log^{-1/2} (f^n)^2}   $ and hence, 
 \begin{align} 
|  \tau ( \psi - \psi_\infty ( \o{f^n \log^{1/2} (f^n)^2}   \  \o{f^n \log^{-1/2} (f^n)^2}   ) ) | 
& \leq  \left( \int_B  \psi_\infty \o{(f^n-f)^2}   \int_B  \psi_\infty \o{ |\nabla(f^n-f)|^2}^{\de,\ka}  
 \right)^{1/2}      
\end{align}  
Hence, 
\begin{align} 
 - A_1  &  \leq - \frac{\mu}{N^4}  + C \frac{ | \nabla u | }{N^4}  \left( \int_B  \psi_\infty \o{(f^n-f)^2}   \int_B  \psi_\infty \o{ |\nabla(f^n-f)|^2}^{\de,\ka}    
 \right)^{1/2}   \\  
     & \leq - \frac{\mu}{N^4}  + C  |\nabla u|^2 \frac{\eta}{N^4} + 
 \frac1{10 {N^4}} \int_B  \psi_\infty \o{
   |\nabla(f^n-f)|^2}^{\de,\ka}.    
\end{align}     
 \vspace{.2cm}

The term between parentheses in the definition of $A_3$ can be written as 
\begin{equation} 
         \overline {  \frac{\nabla  u^n}{f^n}   \Big( \log^{1/2} (f^n)^2 + \log^{-1/2} (f^n)^2   \Big)   
\Big[  f^n \log^{1/2}   (f^n)^2  - \o{  f^n \log^{1/2}   (f^n)^2    }       \Big] }   
\end{equation}     
 
If we denote $\nu_{t,x,R}(\Pi, \lambda)$ the Young measure associated to the 
sequence $ (\nabla_x u^n, f^n) $, then we see easily that $A_3$ is given by 

\begin{align} 
 A_3  & = -    \frac{2ak}{N^4}  \int_B  \int\int 
    \left(   \frac{\Pi}{\lambda} (\log^{1/2} \lambda^2 + \log^{-1/2} \lambda^2    )  -
   \frac{\Pi'}{\lambda'} (\log^{1/2} \lambda'^2 + \log^{-1/2} \lambda'^2    )       \right)    \\ 
 & \quad \quad \quad  \quad \quad      ( \lambda \log^{1/2}\lambda^2  - \lambda' \log^{1/2}\lambda'^2   ) 
  \frac{R_iR_j}{1-|R|^2} \psi_\infty  \, \,   d\nu \, d\nu' \,  dR.       
\end{align}     
The absolute value of the  two factors inside the integral can be bounded respectively by 
\begin{align} 
 & |  {\Pi}   - {\Pi'} | (   \frac{\log^{1/2} \lambda^2}\lambda   +  \frac{\log^{1/2} \lambda'^2}{\lambda'}   ) 
 +  ( |\Pi| + |\Pi'| )  (    \frac{\log^{1/2} \lambda^2}\lambda   -  \frac{\log^{1/2} \lambda'^2}{\lambda'}    )    \quad \hbox{and}     \\ 
& |\lambda -\lambda'| (\log^{1/2} \lambda^2 + \log^{1/2} \lambda'^2  ).  \quad \hbox{Hence}
\end{align}     
   
\begin{align} 
|  A_3 |   & \leq   \frac{1}{10 {N^4}}  \int_B  \int\int 
   |  {\Pi}   - {\Pi'}   |^2 (   \frac{\log  \lambda^2}\lambda   +  \frac{\log  \lambda'^2}{\lambda'}   )^2
\frac{1}{1-|R|^2} \psi_\infty   \,   d\nu \, d\nu' \,  dR   \\ 
& \quad \quad \quad  +  \frac{C}{N^4}  \int_B  \int\int   (|  {\Pi}|   +  |{\Pi'}   | ) |\lambda -\lambda'|^2  
  (   \frac{\log  \lambda^2}{\lambda^2}   +  \frac{\log  \lambda'^2}{\lambda'^2}   )^2  
\frac{1}{1-|R|^2} \psi_\infty   \,   d\nu \, d\nu' \,  dR   \\ 
& \quad \quad  \quad \quad\quad \quad +  \frac{C}{N^4}  \int_B  \int\int   (1 + |  {\Pi}|   +  |{\Pi'}   | )  |\lambda - \lambda'|^2 \frac{1}{1-|R|^2} \psi_\infty   \,   d\nu \, d\nu' \,  dR  \\ 
& \leq  \frac{1}{10 {N^4}}  \mu   +  \frac{1}{10 {N^4}} \ka   +    \frac{C}{{N^4} }  |\nabla u|^2 \eta 
\end{align}   
 \vspace{.2cm}

Finally, to bound $-A_4$, we split it into two terms : 
\begin{align} 
|A_4^1 |  & \leq  \frac2{N^4} \int_B  \psi_\infty  
 | \gamma_{ij}|  | \nabla g | \ dR \\ 
 & \leq   \frac{1}{10 N^4  }  \o{ |\nabla  u^n - \nabla u  |
 | }^{\de,\ka}  + 
  \frac{C}{N^4} 
    \o{ \left( \int_B ( g^n - g  )   |\nabla_R  g| \psi_\infty  \right)^2}     \\
  & \leq   \frac{1}{10 N^4 }  \mu  +    \frac{C}{N^4}  \int_B
    \psi_\infty |\nabla_R g|^2 \int_B \psi_\infty  \o{  ( g^n - g  )^2    }    
\end{align}   
 
To bound $A_4^2   $, we first consider the case $k > 1 $ where the term can be treated 
as $A_4^1 $  using \eqref{Hardy-1}: 
\begin{align} 
\nonumber |A_4^2 |  & \leq  \frac2{N^4} \int_B  \psi_\infty  (| \gamma_{ij} |  + | \gamma_{ij}' |)  \frac{g}{1-|R|}  
    \ dR \\  
 \label{k-big}  & \leq     \frac{1}{10 {N^4}}  \mu  +    \frac{C}{ {N^4} } \int_B \psi_\infty 
    \o{  ( g^n - g  )^2    } \ dR        \int_B \psi_\infty  \frac{ | g|^2 }{ ( 1- |R|)^2} \ dR   \\ 
\nonumber  &  \leq     \frac{1}{10{N^4} }  \mu    +    \frac{C}{N^4}  \left( \int_B \psi_\infty |\nabla_R g|^2 
  \right)   \eta. 
\end{align}

In the case $ k \leq 1 $, we have to
 use  \eqref{Hardy-inter-log} instead of \eqref{Hardy-1}. 
Take $  0\leq \beta < k   $ and $\gamma = \frac{1-\beta}2$

\begin{equation} \label{a42} 
   \left(  \int_B \psi_\infty 
         \frac{ |g^n - g | \,  g   }{ ( 1- |R|^2)} \ dR  \right)^2   
   \quad \leq  \int_B \psi_\infty  \frac{ (g^n -g)^2   \log^{-\gamma} (f^n)^2 }{(1-|R|^2)^{1-\beta}} 
  \, dR \quad    \int_B \psi_\infty  \frac{  g^2   \log^{\gamma} (f^n)^2 }{(1-|R|^2)^{1+\beta}} \ dR  
\end{equation}   
   
To bound the second term, we use the following Young's inequality  for $a,b \geq 1$
$ab \leq a \log^\gamma a + e^{(b^{\frac1\gamma})} $. We  denote $d= 1
- |R|^2$ and hence

\begin{align*} 
  \int_B \psi_\infty  \frac{  g^2   \log^{\gamma} (f^n)^2 }{(1-|R|^2)^{1+\beta}} \ dR   
 &\leq   \int_B \psi_\infty   \left[ \frac{g^2}{d^{1+\beta}}  \log^\gamma  \frac{g^2}{d^{1+\beta}}  
 + |f^n|^2   \right] dR.  
\end{align*}   

On the set $\{  g^2 \geq \frac1{d^\eps}   \}   $    where $\eps = \frac{k-\beta}2$, we have 
$  \log^\gamma  \frac{g^2}{d^{1+\beta}} \leq C \log^\gamma g^2  $. Besides, we have  using 
\eqref{Hardy-inter-log} 

\begin{equation*} 
   \int_B \psi_\infty    \frac{g^2 \log^\gamma  {g^2}  }{d^{1+\beta}}   dR  \ 
\leq   \left( \int_B \psi_\infty   g^2 \log   {g^2}  \right)^{1-\beta \over 2} 
 \left( \int_B \psi_\infty (  |\nabla_R g|^2  + g^2) 
  \right)^{1 + \beta \over 2}.  
\end{equation*}  

On the set $\{  g^2 \leq \frac1{d^\eps}   \}   $, we have 
$$ \frac{g^2}{d^{1+\beta}}  \log^\gamma  \frac{g^2}{d^{1+\beta}}  \leq  \frac{C}{d^{1+\beta+\eps}} 
 \log^\gamma (\frac1d)   $$ 
which is integrable in the ball $B$ with the measure $\psi_\infty dR$.

To bound the first term on the right hand side of \eqref{a42}, we first notice  that 
$$  \o{ (g^n -g)^2   \log^{-\gamma} (f^n)^2 }   \leq 
 \o{  (f^n-f)^2  \log^{1-\gamma} (C +  (f^n -f)^2  )    }     $$
 which can be easily proved using Young measures. Besides, we have  using 
\eqref{Hardy-inter-log}

\begin{align*} 
&   |  \int_B \psi_\infty  \frac{   (f^n-f)^2  \log^{1-\gamma} (C +  (f^n -f)^2  )    }{(1-|R|^2)^{1-\beta}} 
  \, dR | \\
& \quad \leq  \left(  \int_B \psi_\infty  (f^n-f)^2  \log (C +  (f^n -f)^2  )     \right)^{1+\beta \over 2} 
   \left(  \int_B \psi_\infty  |\nabla_R(f^n-f)|^2    \right)^{1-\beta \over 2} 
 \\
& \quad \leq  \frac{C}{\lambda^{2 \over 1+\beta}} 
  \left(  \int_B \psi_\infty  (f^n-f)^2  \log (C +
 (f^n -f)^2  )     \right) +  \lambda^{2 \over 1-\beta} 
   \left(  \int_B \psi_\infty  |\nabla_R(f^n-f)|^2    \right). 
\end{align*}   
for each $\lambda > 0$. Passing to the limit weakly (more precisely, applying
$\o{F_n}^{\de,\ka}$) to both sides and optimizing in 
$\lambda$, we deduce that, 
        
 \begin{equation} \label{a42bis} 
  \frac1{N^4} \o{   \left(  \int_B \psi_\infty 
         \frac{ |g^n - g | \,  g   }{ ( 1- |R|^2)} \ dR  \right)^2 } 
\leq    
  \frac{C}{N^4}  \left( \int_B \psi_\infty   g^2 \log   {g^2}  \right)^{1-\beta \over 2} 
 \left( \int_B \psi_\infty (  |\nabla_R g|^2  + g^2) 
  \right)^{1 + \beta \over 2} \eta^{1+\beta \over 2  } \varpi^{1-\beta
  \over 2}. 
\end{equation}

Putting all these estimates together, we deduce that

  \begin{align}  \label{gron}
(\partial_t + u . \nabla) \frac{ \eta  }{N^4}  + \frac{  \mu +   \varpi   } {4N^4} 
\leq    C |\nabla u|^2 \frac{ \eta  }{N^4}  + \frac{C}{N^4}  \left(1+  \int_B \psi_\infty |\nabla_R g|^2 
  \right)   \Big(  \int_B \psi_\infty  g^2 \log g^2    \Big)^{1-\beta
    \over 1+\beta}  \eta  .     
\end{align}   
       
We can take $\beta = 0$. 
Next, we observe that 
$\int_B \psi_\infty  g^2 \log g^2  dR \leq C  N^2  $. Indeed, if we
introduce 
$h^n = g^n \log^{1/2} g^n  $, we see that $N^n_2 \geq \left( \int_B
\psi_\infty (h^n)^2 \right)^{1/2} $  and then it is easy to see 
using that $(x,y) \to \frac{x^2}y$ is convex that 
$$  \o{\left( \int_B
\psi_\infty (h^n)^2 \right)^{1/2}   }  \geq  \left( \int_B
\psi_\infty h^2 \right)^{1/2}   $$
from which we deduce the claim. Hence \eqref{gron} becomes 

 \begin{align}  \label{gron2}
 \frac{d}{dt} \frac{ \eta   }{N^4} (t,X(t,x) )  + \frac{ \mu  + \varpi
 } {4N^4}
 (t,X(t,x) )  
\leq    C |\nabla u|^2 \frac{ \eta  }{N^4} (t,X(t,x) )   +
  C \left[1+  \int_B \psi_\infty  \frac{ |\nabla_R g|^2 }N  
  \right]  \frac{\eta}{N}(t,X(t,x) )    .     
\end{align}   
First notice that the right hand side of \eqref{gron2} is in 
$L^1((0,T) \times K  ) $ for any bounded measurable set of
$\Omega$ (To prove this, we can observe that $\frac\eta{N}$ is bounded and that 
using \eqref{log2-b-lim}, the term between 
brackets in \eqref{gron2} is in $L^1((0,T) \times K  ) $). Hence 
\eqref{gron2} is well justified  in the sense of distribution.      
In particular this justifies all the calculations  done in
this subsection starting from \eqref{Gron-eq}. 

Now, since  the term between 
brackets in \eqref{gron2} is in $L^1((0,T) \times K  ) $, 
for almost all $x$, $ \int_0^T \left[1+  \int_B \psi_\infty  \frac{ |\nabla_R g|^2 }N  
  \right] (t,X(t,x) )    $ is finite. Besides, for almost all 
$x$, ${N}(t,X(t,x) )  $ (which is constant in $t$) is bounded. 
Hence, we deduce that for almost all $x$,  
$ \int_0^T  N^3 \left[1+  \int_B \psi_\infty  \frac{ |\nabla_R g|^2 }N  
  \right]  + |\nabla u|^2    (t,X(t,x) )    $ is finite.
Hence, by Gronwall lemma, we deduce that for a.e $x$, we have 
for all $t<T$,  
$ \frac{\eta (t,x)}{N^4}  \leq  \frac{\eta (0,x)}{N^4}  e^{C_T(x)}  $
and  since  $ \eta (0,x) = 0  $  due to the initial strong convergence, we deduce that 
$ \frac{\eta (t,x)}N^4 = 0$ and hence  $\eta = 0  $ and 
    we deduce the strong convergence of 
$g^n$ to $g$. This yields that $(u,\psi)$ is a weak solution of \eqref{micro}
with the initial data $(u_0,\psi_0)$.  
 

\section{Approximate  system} \label{app-seq}

In the previous section, we proved the weak compactness of a sequence of 
solutions to the system \eqref{micro}. Of course one has  to construct a 
sequence of (approximate) weak  solutions to which we can apply the strategy 
of the previous sections. The only thing we have to make sure is that 
the calculations  done in the previous section can be made on the 
approximate system.  We consider a sequence of global smooth  solutions 
$(u^n,\psi^n)$ to  the following regularized   system where $k$ is some integer that 
depends on $D$. In particular one can take $k=1$ for $D=2$ or $3$:  

 \begin{equation} \label{app}   \left\{  \begin{array}{l}
  {\partial_t u^n} + (u^n\cdot \nabla) u^n- \nu \Delta u^n + \frac1n (\Delta)^{2k} u^n   + \nabla p^n
 = {\div} \tau^n, \quad {\div} u = 0, \\
\\  
\partial_t \psi^n + u^n. \nabla \psi^n =   {\rm div}_R \Big[ - \nabla u^n  \,  R \psi^n
      + {\beta} \nabla \psi^n +   \nabla \U  \psi^n  \Big]  \\
\\ 
\tau^n_{ij} =   \int_B   (R_i \otimes  \nabla_j \U) 
 \psi^n(t,x,R) dR \,  \quad \quad 
 ( \nabla \U  \psi^n  +  {\beta}  \nabla \psi^n).n = 0  \;  \hbox{on} \;   
 \partial B(0,R_0).       
  \end{array} \right. 
\end{equation}   
with  a smooth 
   initial condition $(u^n_0,\psi^n_0)  $ 
such that 
 $(u^n_0, \psi^n_0)$  converges  strongly to  $(u_0, \psi_0)$ in $L^2(\Omega) \times L^1(\Omega \times B)$ 
 and $ \psi^n_0 \log \frac{\psi^n_0 }{ \rho^n_0  \psi_\infty } -\psi^n_0 + \rho^n_0  \psi_\infty  $ 
  converges strongly to $\psi_0  \log \frac{\psi_0}{\rho_0 
 \psi_\infty}  -\psi_0 + \rho_0 \psi_\infty   $
in $  L^1(\Omega \times B) $. We also assume that 
 \eqref{ini-log2} holds uniformly in $n$. 
  In the case $\Omega$ is a bounded domain of $\R^D$, we also 
impose the following boundary condition $ u^n = \Delta u^n  = ... = (\Delta)^{2k-1} u^n = 0    $
at the boundary $\partial \Omega$.

We do not detail the proof of existence for the   system \eqref{app}. We only mention that
we have to combine classical results about strong solutions to Navier-Stokes system 
with the study of the linear  Fokker-Planck equation (see \cite{Masmoudi08cpam}). 
In particular the following operator was used 
 
\begin{equation}  
L \psi = -div ( \psi_\infty \nabla \frac{\psi}{\psi_\infty}   )
 \end{equation}
on the space $\H = L^2(\frac{dR}{\psi_\infty})$ and with domain 
\begin{equation}  
D(L) = \left\{ \psi \in \H  |  \psi_\infty \nabla \frac{\psi}{\psi_\infty}  \in \H, \quad 
   div ( \psi_\infty \nabla \frac{\psi}{\psi_\infty}   )  \in \H \, \quad \hbox{and} \, 
  \psi_\infty \nabla \frac{\psi}{\psi_\infty}|_{\partial B} = 0   \right\} . 
 \end{equation} 
Also the following two  Hilbert spaces  $\H^1$ and $\H^2$ are used in the construction :  
\begin{eqnarray}   
\H^1 = \left\{\psi \in \H  \   | \quad \int \psi_\infty  \left| \nabla \frac{\psi}{\psi_\infty}\right|^2 
 + \frac{\psi^2}{\psi_\infty}  \ dR  < \infty  \right\}   \\ 
 \H^2 = \left\{\psi \in \H^1 \, 
  | \quad \int  \left(   div ( \psi_\infty \nabla \frac{\psi}{\psi_\infty}  
  )        \right)^2 
   \frac{dR} {\psi_\infty}  < \infty .   \right\} 
 \end{eqnarray} 
Following the the proof of existence given in \cite{Masmoudi08cpam}, we can prove 

\begin{prop} \label{prop-ex}
\label{fene-app}
Take $u_0^n  \in H^s (\Omega) $ and $\psi_0^n \geq 0  $ such that  $ 
 \psi^n_0 - \rho^n_0 \psi_\infty   \in  H^s (\Omega;
  L^2 ( {dR \over \psi_{\infty}}  )  )  $ with  $ \rho^n_0 =  \int \psi_0^n  dR  \in L^\infty (\Omega)$. 
Then, there exists a global   unique solution  
  $(u^n,\psi^n)$ to  (\ref{app}) such that 
$(u^n,\psi^n -\rho^n \psi_\infty )$ is in 
  $C([0,T); H^s) \times C([0,T); H^s (\R^N; L^2 ( {dR \over \psi_{\infty}}  )  ) ) $ for all $0< T$.  Moreover, 
$ u^n \in L^2([0,T);H^{s+k})$ and $\psi^n -\rho^n \psi_\infty  \in  L^2([0,T);H^{s}(\R^N;\H^1) )$.  
\end{prop}

\begin{rem}
The proof is exactly the same as the proof of Theorem 2.1 of \cite{Masmoudi08cpam} with few differences:
\begin{itemize}
\item  In  \cite{Masmoudi08cpam}, we only had local existence  whereas here, 
 we have global existence since we have more regularity. 

\item Theorem 2.1 of \cite{Masmoudi08cpam} 
was stated in the whole space. Of course in the case of a bounded domain, we have to use 
energy bounds for Navier-Stokes written in a bounded domain. 

\item In theorem 2.1 of \cite{Masmoudi08cpam}  we assumed that $\int \psi_0 dR = 1$. 
The result can be easily extended to this  more general case.  We also point out that there is 
a small mistake in  the 
statement of the theorem 2.1 of \cite{Masmoudi08cpam}. Indeed, one has to read
$\psi_0 - \psi_\infty \in H^s (\Omega;
  L^2 ( {dR \over \psi_{\infty}}  )  )  $
instead of $\psi_0  \in H^s (\Omega;
  L^2 ( {dR \over \psi_{\infty}}  )  )  $ when the problem is in the whole space. 
\end{itemize}
\end{rem} 
 
It is   clear that the  solutions constructed in Proposition \ref{prop-ex} 
satisfy the free-energy bound \eqref{free} and the 
extra bound \eqref{log2-b}  (with $\Omega$ replaced by $K$ in the whole space case).  

\vspace{.2cm}
 Once we have our sequence of  regular approximate solutions, we have to check 
that all the computations performed in the previous section can be done 
on this sequence $(u^n,\psi^n)$.  The only point to be checked is that 
Proposition \ref{mu-beta} still holds since the rest of the proof only involves 
the transport equation.  Now, 
$v^n$ and $w^n $solve 
\begin{align} \label{v-n-1}
 \left\{     \begin{array}{l} 
  \partial_t v^n - \Delta v^n   + \frac1n \Delta^{2k}  v^n + \nabla p_1^n = \nabla. \tau^n  \\ 
   v^n (t=0) = 0  
  \end{array}\right. 
\end{align}   
\begin{align}\label{w-n-1}
 \left\{     \begin{array}{l} 
  \partial_t w^n - \Delta w^n     + \frac1n \Delta^{2k}  w^n      + \nabla p_2^n =  - u^n.\nabla u^n   \\ 
   v^n (t=0) = u^n(t=0) 
  \end{array}\right. 
\end{align}   
and we define $   v^{n,\De}$ the solution of 
\begin{align} \label{v-n-d-k-1}
 \left\{     \begin{array}{l} 
  \partial_t v^{n,\De} - \Delta v^{n,\De} +    \frac1n \Delta^{2k}   v^{n,\De} +     \nabla p_1^{n,\De} = \nabla. \tau^{n,\De}  \\ 
   v^{n,\De} (t=0) = 0   
  \end{array}\right.  
\end{align}  
Step 1 of the proof of Proposition \ref{mu-beta} is the same with the difference that 
one has  to apply  parabolic regularity for the perturbed Stokes operator which 
yields the same uniform in $n$ estimate. 
Hence, we deduce that 
$ \| \nabla v^{n,\De} -  \nabla v^{n}       \|_{L^p( (0,T) \times \Omega    )  }  $ goes to zero when 
$\De$ goes to zero uniformly in $n$ for $p< 2$.

For the second step, we first notice that \eqref{v-n-d-k-limit} remains the same since 
$\frac1n (\Delta)^{2k} u^n $ converges weakly to zero. Moreover, multiplying 
the first equation of \eqref{v-n-1} by $v^{n,\De}$, we get 

\begin{equation} \label{v-n-d-k-limit2-app-n}
   \partial_t \frac{|v^{n,\De}|^2}2  - \Delta \frac{|v^{n,\De}|^2}2  +  
|\nabla v^{n,\De} |^2  +    G^{n,\De}    + \div ( p_1^{n,\De} v^{n,\De} )   = 
 \div ( v^{n,\De} .  \tau^{n,\De}) -   \tau^{n,\De} : \nabla v^{n,\De} . 
\end{equation}  
 where $G^{n,\De}$ is given by 

\begin{equation}
\begin{array}{l} 
G^{n,\De} = \frac1n [ \div_i( \nabla_i \Delta^{2k-1}  v^{n,\De} . v^{n,\De}  - 
  \Delta^{2k-1}  v^{n,\De} . \nabla_i   v^{n,\De}  +  \\ 
\quad \quad \quad \quad \quad \quad   \nabla_i  \Delta^{2k-2}  v^{n,\De} . 
   \Delta  v^{n,\De}  -...  -  \Delta^{k}  v^{n,\De} . \nabla_i    \Delta^{k-1}  v^{n,\De}      )
 + \Delta^k    v^{n,\De} .\Delta^k    v^{n,\De}  
    ] 
\end{array}
\end{equation}  
Using the fact that $\frac1n \int_0^T \int_\Omega |\Delta^k    v^{n,\De}  |^2  $ and 
$  \int_0^T \int_\Omega |\nabla     v^{n,\De}  |^2  $  are 
uniformly bounded, we deduce easily that $\o{G^{n,\De}} = \o{  |\Delta^k    v^{n,\De}  |^2    } \geq 0  $
and hence passing to the limit in \eqref{v-n-d-k-limit2-app}, we deduce that  

\begin{equation} \label{v-n-d-k-limit2-app}
   \partial_t \frac{|v^{\De}|^2}2  - \Delta \frac{|v^{\De}|^2}2  +  
|\nabla v^{\De} |^2  + \mu_\De    + \div ( p_1^{\De} v^{\De} )   \leq 
 \div ( v^\De .  \tau^{\De}) -   W^{\De \De}. 
\end{equation}  
and hence Proposition \ref{mu-beta} is replaced by an inequality 
 $  \mu  \leq  -  \int_B \beta_{ij}  \frac{R_i R_j}{1-|R^2|} dR $ 
 which is the inequality that we need in the rest of the proof.

\section{Conclusion} \label{conc}
In this paper we gave a proof of existence of weak solutions to the system 
\eqref{micro}, using the fact that a sequence of regular solutions to the 
approximate system \eqref{app} converges weakly to a weak solution  of \eqref{micro}. 
We would like here to mention few important open problems  (with  increasing level of 
difficulty,  at least  this is what  the author thinks): 
 
\begin{itemize} 
 \item {\it The zero diffusion limit in $x$.} If we add a diffusion term 
 $\frac1n \Delta_x \psi$ in the Fokker-Planck equation of \eqref{micro}, then 
 one can prove the global existence of weak solutions to the regularized model. 
 A natural question is whether we recover a weak solution of the unregularized 
 system  \eqref{micro} when $n$ goes to zero. This is the object of a forthcoming paper 
 \cite{Masmoudi10prep-x}. The difficulty comes from the fact that  the calculation of 
section \ref{weak} used in a critical way the fact that we had a transport equation in 
the $x$ variable.

\item {\it Relaxing the assumption \eqref{log2-b}.} This extra bound was only used 
to give some extra control on the stress tensor. Can we prove the same existence 
result without it ? 

\item {\it Other  models.} A natural question is whether we can extend this to the 
Hooke model (where the system can be reduced to a macroscopic model). We were not 
able to perform this. The main difficulty is that we do not know whether 
the extra stress tensor $\tau$ is in $L^2$.  Nevertheless, we know how to use 
the strategy to this paper to prove global existence for the FENE-P model 
\cite{Masmoudi10prep-p}.

\item {\it Regularity in 2D.} Many works on polymeric flows are motivated by similar 
known results for the Navier-Stokes system. In particular a natural question is 
whether one  can prove global existence of smooth solutions to \eqref{micro} in 2D. 
We point out that this is known for the co-rotational model \cite{LZZ08,Masmoudi08cpam}. 
Of course this seems to be a very difficult problem since, we only have an $L^2$ bound 
on $\tau$ and that an $L^\infty$ bound on $\tau$ was necessary in the previously mentioned 
works. In particular the similar result is not known for the co-rotational Oldroyd-B
model where one can  prove $L^p$ bounds on $\tau$  for each $p > 1$.

\item {\it Is system \eqref{micro} better behaved  than Navier-Stokes.} One  
does not expect   to prove results on \eqref{micro} which are not known for Navier-Stokes 
since \eqref{micro} is more complicated than Navier-Stokes. However, one can 
speculate that due to the polymers and  the extra stress tensor, system \eqref{micro} may behave  
better than  Navier-Stokes  and that one can prove global existence of smooth solutions 
to \eqref{micro}  even if such result is not proved or disproved for the Navier-Stokes 
system.

\end{itemize}

\section{Acknowledgments}
The work of N. M.  is partially supported by NSF-DMS grant  0403983.
The author would like to thank P.-L. Lions and Ping Zhang for many 
discussions about this model. He also would like to  thank the IMA 
where part of this work was done.

\end{document}